\documentclass[11pt,reqno]{amsart}

\usepackage{ebgaramond}
\usepackage{hyperref}
\hypersetup{hypertex=true, colorlinks=true, linkcolor=black, anchorcolor=black, citecolor=black}

\usepackage[capitalize, noabbrev, nameinlink]{cleveref} 

\usepackage{amsthm}
\usepackage{amsmath}
\usepackage{amssymb}
\usepackage{mathpazo}
\usepackage{mathrsfs}
\usepackage{mathtools}
\usepackage{libertine}

\usepackage{enumerate}
\usepackage[table, svgnames]{xcolor} 

\usepackage{verbatim}
\usepackage[all]{xy} 
\usepackage{tikz} 
\usepackage{tikz-cd}
\usetikzlibrary{shapes,arrows}

\newtheorem{theorem}{Theorem}[section]
\newtheorem{lemma}[theorem]{Lemma}
\newtheorem{corollary}[theorem]{Corollary}
\newtheorem{proposition}[theorem]{Proposition}
\newtheorem{claim}[theorem]{Claim}

\theoremstyle{definition}
\newtheorem{definition}[theorem]{Definition}

\newcommand\bp{{\bar\partial}}

\renewcommand{\bar}{\overline}
\newcommand{\eps}{\varepsilon}

\newcommand{\beq}{\begin{equation}}
\newcommand{\eeq}{\end{equation}}
\newcommand{\ee}{\end{eqnarray*}}
\newcommand{\be}{\begin{eqnarray*}}

\newcommand{\bd}{\begin{enumerate}}
	\newcommand{\ed}{\end{enumerate}}

\renewcommand{\hat}{\widehat}
\renewcommand{\tilde}{\widetilde}

\renewcommand{\bp}{\bar{\partial}}

\newcommand{\QQ}{{\mathbb Q}}
\newcommand{\RR}{{\mathbb R}}

\newcommand\cB{{\mathcal B}}

\newcommand\cE{{\mathcal E}}

\newcommand\cH{{\mathcal H}}

\newcommand\cO{{\mathcal O}}

\newcommand{\Psh}{\mathrm{Psh}}

\newcommand{\Ent}{\mathrm{Ent}}

\newcommand{\tr}{\mathrm{tr}}

\newcommand{\Ric}{\mathrm{Ric}}
\renewcommand{\d}{\mathrm{d}}

\newcommand{\reg}{\mathrm{reg}}

\newcommand{\Aut}{\mathrm{Aut}}

\newcommand{\PSH}{\mathrm{PSH}}

\newcommand{\ddc}{\mathrm{dd}^{\rm c}}

\DeclareFontFamily{U}{mathx}{}
\DeclareFontShape{U}{mathx}{m}{n}{<-> mathx10}{}
\DeclareSymbolFont{mathx}{U}{mathx}{m}{n}
\DeclareMathAccent{\widehat}{0}{mathx}{"70}
\DeclareMathAccent{\widecheck}{0}{mathx}{"71}

\allowdisplaybreaks[4]

\makeatletter
\newcommand{\constant}[2]{%
  \@ifundefined{c@#1counter}{\newcounter{#1counter}}{}%
  \refstepcounter{#1counter}%
  \protected@write\@auxout{}{%
    \string\newlabel{#2}{{\arabic{#1counter}}{\thepage}{}{#1@\arabic{#1counter}}{}}%
  }%
  #1_{\arabic{#1counter}}%
}
\makeatother

\usepackage[toc,page]{appendix}

\makeatletter
\renewcommand{\tocsection}[3]{%
  \indentlabel{\@ifnotempty{#2}{\bfseries\ignorespaces#1 #2\quad}}\bfseries#3}
\renewcommand{\tocsubsection}[3]{%
  \indentlabel{\@ifnotempty{#2}{\ignorespaces#1 #2\quad}}#3}

\newcommand\@dotsep{4.5}
\def\@tocline#1#2#3#4#5#6#7{\relax
  \ifnum #1>\c@tocdepth 
  \else
    \par \addpenalty\@secpenalty\addvspace{#2}%
    \begingroup \hyphenpenalty\@M
    \@ifempty{#4}{%
      \@tempdima\csname r@tocindent\number#1\endcsname\relax
    }{%
      \@tempdima#4\relax
    }%
    \parindent\z@ \leftskip#3\relax \advance\leftskip\@tempdima\relax
    \rightskip\@pnumwidth plus1em \parfillskip-\@pnumwidth
    #5\leavevmode\hskip-\@tempdima{#6}\nobreak
    \leaders\hbox{$\m@th\mkern \@dotsep mu\hbox{.}\mkern \@dotsep mu$}\hfill
    \nobreak
    \hbox to\@pnumwidth{\@tocpagenum{\ifnum#1=1\bfseries\fi#7}}\par 
    \nobreak
    \endgroup
  \fi}
\def\l@subsection{\@tocline{2}{0pt}{2.5pc}{5pc}{}}
\def\l@subsubsection{\@tocline{3}{0pt}{5pc}{5pc}{}}
\makeatother

\usepackage{fancyhdr}
\begin{document}

\title{Existence of singular cscK metrics on klt varieties}


\author{Yaxiong Liu}
\address{Mathematical Science Research Center, Chongqing University of Technology, No. 69, Hongguang Avenue, Banan District, Chongqing 400054, People's Republic of China}
\email{yxliu269@cqut.edu.cn}
\date{}

\begin{abstract} 

On a polarized klt variety $(X,L)$ with $G=\mathrm{Aut}_0(X)$  the identity component of the automorphism group of $X$, suppose Mabuchi functional is $G$-coercive in $c_1(L)$,
then there exists a singular cscK metric in $c_1(L)$, which is a minimizer of Mabuchi functional.
\end{abstract}
\maketitle

\tableofcontents

\section{Introduction}

Finding a canonical metric on K\"ahler manifolds is the classical, central problem in K\"ahler geometry, arising from Yau's breakthrough of Calabi conjecture \cite{Yau78a}.
In particular, the interest one is the constant scalar curvature K\"ahler (cscK) metric.
Tian's properness conjecture \cite{TG00, Tian94}  predicts that the existence of cscK metric is equivalent the properness of Mabuchi functional.
The conjecture has been solved by \cite{DR17} for K\"ahler--Einstein metrics, \cite{CC21b}, \cite{BDL20} for cscK metrics (see also \cite{AJL23}, \cite{HLiu23}, \cite{DNJL26} for the general weighted cscK metrics).
It provides an analytic and variational characterization of the existence of cscK metrics.

In the study of Minimal Model Program and moduli problem (e.g, of general type or Fano type vareities), it is inevitable to solve canonical metric equations on varieties with singularities (\cite{EGZ09, SSY16, BBEGZ19, LWX19, LTW21, LTW22, Li22b, DNGG23, CDL26} and so on).
Therefore, we are also interested in finding cscK metrics on projective varieties with klt singularities.

First attempt for singular cscK metrics  is due to \cite{PTT23}, in which the authors proved that the coercivity of Mabuchi functional implies the existence of singular cscK metrics on a $\QQ$-Gorenstein smoothable K\"ahler variety with log terminal singularities.
Boucksom--Jonsson--Trusiani \cite{BJT26} established the openness of coercivity of (weighted) Mabuchi functional on the resolution under the appropriate assumption of resolution, the so-called \textit{resolution of Fano type}, which assumption arose from the RCD theory of singular K\"ahler-Einstein variety \cite{Sze25}.  
Under the same assumption with \cite{BJT26}, Pan-T\^{o} \cite{PT26} established the existence of (weighted) cscK metrics provided the coercivity of (weighted) Mabuchi functional.

\medskip

The goal of the article is to establish the existence of singular cscK metrics on the general klt variety without any assumption of singularities.
For now, we briefly state our setting.

Let $X$ be a compact connected klt K\"ahler variety of complex dimension $n\geq2$ with an ample line bundle $L$, and let $\omega\in c_1(L)$ be a smooth K\"ahler form (as a restriction
of a smooth K\"ahler form defined on an local ambient space). 
Let $G=\Aut_0(X)$ be the identity component
of the automorphism group of $X$.
The singular cscK metric in $[\omega]$ on $X$ is defined in Definition~\ref{def:cscK}.
The Mabuchi functional $M$ in $[\omega]$ is defined by \eqref{eq:def-Mabuchi}. 
We define the $G$-coercivity of $M$ in Definition~\ref{def:G-coercive}.

Our main theorem is the following

\begin{theorem}
\label{thm:main}
	
	Suppose Mabuchi functional $M$ is $G$-coercive, then there exists a cscK metric in $c_1(L)$, which also minimizes $M$ in $[\omega]$.
\end{theorem}

We establish the existence via a detailed approximation scheme on a log resolution of $(X,-K_X)$, adapted from \cite[Section~4.1]{CDL26} and  \cite[Section~4.3]{GSS26}. After resolving the singularities, we perturb the pullback of the parabolic system and derive a priori estimates that are uniform in the perturbation parameter (see Section~\ref{sec:Uniform esitmates for perturbed equations}).
A crucial point is that the Mabuchi functional of perturbed K\"ahler class on the resolution may not be $G$-coercive, even though $M$ is $G$-coercive.
To produce the solution of perturbed equations, instead, we use Schauder fixed point theorem by constructing a continuous, relatively compact map on the space $C^{2,\gamma}$, whose fixed point is exactly the solution of the perturbed equation. 

\medskip

The rest of the paper is organized as follows. In Section~\ref{sec:Preliminaries}, we review the geometric pluripotential theory on klt varieties.
In Section~\ref{sec:Uniform esitmates for perturbed equations}, we develop the perturbation equation and establish uniform estimates to guarantee that the limit of perturbed solution is the singular cscK.
In Section~\ref{sec:Proof of the main theorem}, we prove the existence of perturbed equation \eqref{eq:eps_twist_cscK} by employing Schauder fixed point theorem.
Then we prove the Theorem~\ref{thm:main} by combining the uniform estimates and existence of perturbed equation.  

\medskip

\noindent \textbf{Acknowledgments.} 
The author was partially supported by NSFC grant 12601149.
Part of this work was carried out while the third-named author was visiting the Institute for Theoretical Sciences at Westlake University, which he would like to thank for the hospitality and support.

\medskip
\noindent \textbf{Use of AI.} 
This project was Motivated by author's previous work \cite{HLiu25}. 
Adapted from the perturbation method in \cite[Section~4.1]{CDL26} and  \cite[Section~4.3]{GSS26}, 
we developed the perturbation scheme presented in Section~\ref{sec:Uniform esitmates for perturbed equations} and subsequently established these apriori estimates in Section~\ref{sec:Uniform esitmates for perturbed equations}. 

The difficulty is how to prove the existence of the solution of the perturbed system \eqref{eq:eps_twist_cscK} without any assumption of singularity of $X$.
During discussions between the authors and GPT-6 Astra concerning this difficulty,
GPT-6 Astra suggested the Schauder fixed point theorem to solve the equation \eqref{eq:eps_twist_cscK}.
Theorem~\ref{thm:exist-perturb-equ} was critically assisted by GPT-6 Astra.

All mathematical arguments were written by the authors, who take full responsibility
for the content of this paper.

\section{Preliminaries}
\label{sec:Preliminaries}

We recall some basic notions and results for the pluripotential theory and variational theory of cscK metrics on klt varieties.

\subsection{Potential theory on klt K\"ahler space}

Let $X$ be a compact connected klt K\"ahler variety of complex dimension $n\geq2$ with an ample line bundle $L$, and let $\omega\in c_1(L)$ be a smooth K\"ahler form (as a restriction
of a smooth K\"ahler form defined on an local ambient space). 
Let $\mathcal H_\omega$ denote the corresponding space of K\"ahler potentials, using the convention $\ddc = \frac{i}{2\pi} \partial \bar \partial$:
$$\mathcal H_\omega \coloneqq \{u \in C^\infty(X) \colon \omega_u \coloneqq \omega + \ddc u > 0 \}.$$
We denote the normalized K\"ahler potentials
\begin{align*}
	\cH_0\coloneqq\{u\in\cH_\omega\colon \sup_Xu=0\}.
\end{align*}
By $\PSH_\omega$ we denote the corresponding space of $\omega$-psh functions.
We refer to \cite{Dem85} for more details.

Following \cite{BEGZ10,BBEGZ19}, for any $u\in\PSH_\omega$, the
associated non-pluripolar Monge--Amp\`ere measure is defined by
\begin{align*}
    \omega_u^n
    \coloneqq
    \lim_{k\to\infty}
    \mathbf{1}_{\{u>-k\}}
    \left(\omega+\ddc\max\{u,-k\}\right)^n.
\end{align*}
This is a positive Borel measure on $X$ satisfying
\[
    \int_X\omega_u^n\leq V_\omega
    \coloneqq\int_X\omega^n.
\]
For $p\geq1$, the full-mass class and the finite-energy class are
defined by
\begin{align*}
    \cE_\omega\coloneqq
    \left\{
        u\in\PSH_\omega
        \ \middle|\
        \int_X\omega_u^n=V_\omega
    \right\}, \ \ \  \ \ 
    \cE^p_\omega\coloneqq
    \left\{
        u\in\cE_\omega
        \ \middle|\
        \int_X|u|^p\omega_u^n<\infty
    \right\}.
\end{align*}
Thus we have the strict inclusions 
\[
    \cH_\omega \subset\cE^p_\omega \subset\cE_\omega \subset\PSH_\omega.
\]
It follows from \cite{DNG18} that, for every $p\geq1$, the $L^p$ type Mabuchi metric $d_p$ on $\cE^p_\omega$ is complete, moreover 
$\cH_\omega \subset \cE^p_\omega$ is $d_p$-dense. Additionally, $\bigl(\cE^p_\omega,d_p\bigr)$ is a
geodesic metric space, with a distinguished class of geodesic segments coming from solutions to a degenerate complex Monge--Amp\`ere equation.  The case $p=1$ was obtained independently in
\cite{Dar17a}. 

For $u\in\cH_\omega$, the Monge--Amp\`ere energy is defined by
\[
    E(u) \coloneqq \frac{1}{(n+1)V_\omega}\sum_{j=0}^n \int_X u\,\omega_u^j\wedge\omega^{n-j}.
\]
The same formula extends to bounded $\omega$-psh potentials. Using the
monotonicity of $E$, one then defines, for any $v\in\PSH_\omega$,
\[
    E(v) \coloneqq \inf\left\{ E(u) \ \middle|\ u\in\PSH_\omega\cap L^\infty(X),\ u\geq v \right\} \in[-\infty,\infty).
\]
We have $E(v)>-\infty$ if and only if
$v\in\cE^1_\omega$. Moreover, $E$ is affine along finite-energy weak
geodesics, reflecting the close relation of the Monge--Amp\`ere energy to the $d_1$-geometry of $\cE^1_\omega$. \medskip

\subsection{The adapted measure and cscK metrics.}

Let $r$ be a positive integer such that $rK_X$ is Cartier. Let $h$ be a smooth Hermitian metric on the line bundle $rK_X$ such that. For a nowhere vanishing holomorphic section $\Omega$ of $rK_X$ over an open set $U$,  $|\Omega|_{h}$ is a smooth function on $U$.
The adapted measure is locally defined by
\begin{align*}
    \mu_X\coloneqq \frac{\big(i^{rn}\Omega\wedge\bar{\Omega}\big)^{\frac{1}{r}}}{|\Omega|_{h}^{\frac{2}{r}}}.
\end{align*}
We may assume that $\int_X\mu_X=1$.
The  klt condition of $X$ yields that the adapted measure $\mu$ has $L^{p_0}$-density relative to $\omega^n$ for some $p_0>1$ (see \cite[Lemma~3.2]{BBEGZ19}).
More precise, given a smooth resolution $\pi:Y\rightarrow X$, 
then
\begin{align*}
	\pi^*\mu_X=C\prod_i|s_i|_{h_i}^{2a_i}\d V,\qquad
	a_i>-1,
\end{align*}
where $\d V$ is a smooth positive volume form, $s_i$ is a defining section of the exceptional divisor $E_i$, and $h_i$ is a smooth Hermitian metric on $\cO(E_i)$.
Its Ricci form is defined as
\begin{align*}
	\Ric(\mu_X)\coloneqq\frac{1}{r}\ddc\log|\Omega|_h^2,
\end{align*}
which is a smooth form on $X$. There exists $K_0>0$ such that
\begin{align}
\label{eq:Ric(mu_X)-bound}
	-K_0\omega\leq \Ric(\mu_X)\leq K_0\omega.
\end{align}
For a smooth K\'ahler metric $\omega_0$ on $X$, its Ricci curvature is given by
\begin{align*}
	\Ric(\omega_0)\coloneqq -\ddc\log\left(\frac{\omega_0^n}{\mu_X} \right)+\Ric(\mu_X).
\end{align*}

\begin{definition}
\label{def:cscK}
	We say that $\omega_{u}$ is a a \textit{singular cscK} metric if $\omega_{u}$  is a genuine K\"ahler metric on $X_{\reg}$ and $u\in\PSH_\omega\cap C^0(X)$ satisfies the following equation on $X_\reg$
\begin{align}
\label{eq:def of singular cscK}
	\bar s \omega_{u}^n=n\Ric(\omega_{u})\wedge\omega_{u}^{n-1},
\end{align}
where $\Ric(\omega_{u})=-\ddc\log\frac{\omega_{u}^n}{\mu_X}+\Ric(\mu_X)$ and $\bar s=\frac{n}{V_\omega}\int_X\Ric(\mu_X)\wedge\omega_X^{n-1}$.
\end{definition}

For any $u\in\PSH_\omega\cap L^\infty$, the Mabuchi functional is defined by
\begin{align}
\label{eq:def-Mabuchi}
	M(u)\coloneqq \Ent_{\mu_X}(V_\omega^{-1}\omega_u^n) -E_{\Ric(\mu_X)}(u)+\bar{s}E(u),
\end{align}
where
\begin{align*}
	\Ent_{\mu_X}(V_\omega^{-1}\omega_u^n)
	\coloneqq V_\omega^{-1}\int_X\log\left(\frac{V_\omega^{-1}\omega_u^n}{\mu_X}\right)\omega_u^n
\end{align*}
and
\begin{align*}
	E_{\Ric(\mu_X)}(u)\coloneqq
	\frac{1}{V_\omega}\sum_{k=0}^{n-1}\int_Xu\Ric(\mu_X)\wedge\omega_u^k\wedge\omega^{n-1-k}.
\end{align*}

Let $G=\Aut_0(X)$ be the identity component of automorphism group of $X$.
By \cite[Lemma~3.53]{BJT26}, for any $\sigma\in G$ and any $\omega_\varphi\in [\omega]$, then there exists $\varphi_\sigma\in C^\infty(X)$ such that $\sigma^*(\omega_\varphi)=\omega_\varphi+\ddc\varphi_\sigma$.
We denote $\sigma\cdot\varphi=\varphi^\sigma=\varphi+\varphi_\sigma$.the action on K\"ahler potential, that is, $\omega+\ddc\varphi^\sigma=\sigma^*(\omega+\ddc\varphi)$.
We refer to \cite[Section~3.8]{BJT26} for a detailed discussion on the action.

One defines the reduced $\d_1$-distance
\begin{align*}
	\d_1(\phi,\varphi)\coloneqq
	\inf_{\sigma\in G}\d_1(\phi,\sigma^\sigma).
\end{align*}

\begin{definition}
\label{def:G-coercive}
	We say that the Mabuchi functional $M$ is $G$-\textit{coercive} if there exist constant $\lambda_c, K_c>0$ such that for any $v\in\cH_0$
	\begin{align}
	\label{eq:coercive-M}
		M(v)\geq \lambda_c \d_{1,G}(0,v)-K_c.
	\end{align}
\end{definition}

\begin{lemma}
\label{lem:G-inv-Mabuchi}
	If $M$ is bounded from below, then $M$ is $G$-invariance.
\end{lemma}
\begin{proof}
	By \cite[Section~3.8.6]{BJT26}, $M$ is quasi-invariance, that is, for any $\varphi\in\cH$ and any $\sigma\in G$, there exists $\xi(\sigma)\in\RR$ such that
	\begin{align*}
		M(\varphi^\sigma)-M(\varphi)=\xi(\sigma).
	\end{align*}
	Then for any $k\geq 1$,
	\begin{align*}
		\xi(\sigma^k)
		=M(0^{\sigma^k})-M(0)
		=\sum_{j=0}^{k-1}\Big(M(0^{\sigma^{j+1}})-M(0^{\sigma^j}) \Big)
		=k\xi(\sigma)
		\end{align*}
		Note that $\xi(\sigma\circ\sigma')=\xi(\sigma)+\xi(\sigma')$. Then $\xi(\sigma^{-1})=-\xi(\sigma)$.
		
		Denote $K_{\min}=\inf_\cH M>-\infty$. Then for any $k\geq1$,
		\begin{align*}
			K_{\min}\leq M(0^{\sigma^k})
			=M(0)+k\xi(g),\quad
			K_{\min}\leq M(0^{\sigma^{-k}})
			=M(0)-k\xi(g).
		\end{align*}
		Let $k\rightarrow\infty$, we obtain $\xi(\sigma)=0$ for any $\sigma\in G$. We done.
\end{proof}

\section{Uniform esitmates for perturbed equations}
\label{sec:Uniform esitmates for perturbed equations}
In this section, we develop a detailed approximation scheme on a log resolution of $(X,-K_X)$. Then we perturb the pullback of cscK equation and establish the uniform estimates.

\bigskip

The equation \eqref{eq:def of singular cscK} can be re-written as a couple system
\begin{equation}
\label{eq:cscK_system on X}
\begin{cases}
     \omega_{u}^n
		= \; e^{F}  \mu_X, \\
  \Delta_{u} F
		=-\bar s + \tr_u\Ric(\mu_X). 
\end{cases}
\end{equation}\smallskip

Given a smooth log resolution $\pi:Y\rightarrow X$ with exceptional divisor $E \subset Y$.
Slightly abusing precision, we will often not distinguish $X_\reg$ and $Y\backslash E$, and we will often identify the spaces $\mathcal{H}_\omega(X)/\mathcal{E}^1_\omega(X)/\PSH_\omega(X)$ with $\mathcal{H}_{\pi^* \omega}(Y)/\mathcal{E}^1_{\pi^* \omega}(Y)/\PSH_{\pi^* \omega}(Y)$. 

Take the pullback of the system \eqref{eq:cscK_system on X}, we obtain 
\begin{equation}
\label{eq:cscK_system on Y}
\begin{cases}
     (\pi^*\omega+\ddc u)^n
		= e^{\pi^*F} \pi^*\mu_X, \\
  \Delta_{u} (\pi^*F)
		=-\bar s + \tr_{u}\pi^*\Ric(\mu_X).
\end{cases}
\end{equation}\smallskip
The equation \eqref{eq:cscK_system on Y} admits the degenerate class $[\pi^*\omega]$.
The strategy of solving the equation \eqref{eq:cscK_system on Y} 
is to  perturb the semiample class $[\pi^* \omega]$ to a K\"ahler class, and at the same time perturb the singular measure $\pi^*\mu_X$  to become smooth. Then we will estibilish apriori estimates for the resulting perturbed system that are preserved in the limit.  
\medskip

We fix a small K\"ahler form $\omega_Y$.
For simplicity, we always denote $\hat\omega\coloneqq\pi^*\omega$ and $\theta\coloneqq\pi^*\Ric(\mu_X)$.
We introduce the perturbed background K\"ahler metric:
\begin{align*}
	\omega_\eps\coloneqq\hat\omega+\eps\omega_Y.
\end{align*} 
We may assume that $\omega_\eps\leq \kappa\omega_Y$ for some $\kappa>0$ and $0<\eps\ll1$.

The adapted measure $\pi^*\mu_X$ can be regularized by the smooth measures
\begin{align}
	\mu_\delta\coloneqq
	C_\delta e^{h_0}\prod_i(|s_i|_{h_i}^2+\delta)^{a_i}\omega_Y^n,\qquad
	\int_Y\mu_\delta=1,
\end{align}
where $h_0=\log(\frac{\d V}{\omega_Y^n})$. We denote $f_\delta\coloneqq \frac{\mu_\delta}{\omega_Y^n}$.
Then
\begin{align}
\label{eq:Lp-int-mu-delta}
	\left\|f_\delta\right\|_{L^{p_0}(\omega_Y)}\leq C.
\end{align}

Adapted from \cite[Section~4.1]{CDL26} (also \cite[Section~4.3]{GSS26}) , we consider the perturbed twisted cscK equation
 \begin{align}
 \label{eq:eps_twist_cscK} 
 	\begin{cases}
V_{\omega_\eps}^{-1}(\omega_\eps+\ddc{u}_{\eps,\delta,t})^n
 	=e^{F_{\eps,\delta,t}}\mu_\delta,
  \\[1mm]
\omega'\coloneqq \omega_\eps+\ddc{u}_{\eps,\delta,t}, \quad  \sup_Yu_{\eps,\delta,t}=0, \\[1mm]
\Delta_{\omega'}F_{\eps,\delta,t}
=\tr_{\omega'}\theta_{\eps,t}-c_{\eps,t},
\end{cases}
 \end{align}
 where $\theta_{\eps,t}\coloneqq t\theta-(1-t)A\omega_\eps$ for some $A>0$, $t\in[0,1]$ and
 \begin{align*}
 	c_{\eps,t}\coloneqq \frac{n}{V_\eps}\int_Y\theta_{\eps,t}\wedge\omega_\eps^{n-1}.
 \end{align*}
 By \eqref{eq:Ric(mu_X)-bound}, there exists a constant $\constant{K}{constK:theta_eps,t-bound}>0$ independent of $\eps,t$, such that
 \begin{align}
 \label{eq:theta_eps,t-bound}
 	-K_{\ref{constK:theta_eps,t-bound}}\omega_\eps\leq \theta_{\eps,t}
 	\leq K_{\ref{constK:theta_eps,t-bound}}\omega_\eps, \quad
 	|c_{\eps,t}|\leq nK_{\ref{constK:theta_eps,t-bound}}.
 \end{align}
 Indeed, the upper bound is due to \eqref{eq:Ric(mu_X)-bound} and the definition of $\theta_{_\eps,t}$.
 For the lower bound, again by \eqref{eq:Ric(mu_X)-bound}, one has
 \begin{align*}
 	-t\pi^*\Ric(\mu_X)\leq tK_0\hat\omega\leq K_0\hat\omega,
 \end{align*}
 which yields 
 \begin{align*}
 	-\theta_{\eps,t}=-t\pi^*\Ric(\mu_X)+(1-t)A\omega_\eps
 	\leq (K_0+A)\omega_\eps.
 \end{align*}
 
For ease of notation, a constant is said to be \textit{uniform} if it is independent of $\eps$, $\delta,t$. We always use the constant $K$ to denote a uniform constant. Sometimes, we use subscript to emphasize its dependence (e.g. $K_{\theta}$ depends on $\theta$).

\subsection{Uniform $C^0$-estimate}
In this subsection, we establish the uniform $C^0$-estimate of the perturbed equation \eqref{eq:eps_twist_cscK} provided the uniform entropy bound.

We denote
\begin{align*}
	H_{\eps,\mu_\delta}(u)
	=\Ent_{\mu_\delta}(V_\eps^{-1}\omega_{\eps,u}^n).
\end{align*}

\begin{theorem}
\label{thm:unif-C0}
	Suppose $(u_{\eps,\delta,t}, F_{\eps,\delta,t})$ is a solution of \eqref{eq:eps_twist_cscK} and $H_{\eps,\mu_\delta}(u_{\eps,\delta,t})=\int_YF_{\eps,\delta,t}e^{F_{\eps,\delta,t}}\mu_\delta\leq H$.
	Then there exists a constant $K_H=K(H,\omega_Y)$ such that
	\begin{align*}
		\|u_{\eps,\delta,t}\|_{L^\infty}+\| F_{\eps,\delta,t}\|_{L^\infty}\leq K_H.
	\end{align*}
\end{theorem}
	For simplicity, we omit the subscript $u=u_{\eps,\delta,t}, F=F_{\eps,\delta,t}$.

\begin{lemma}
	There exists a constant $\constant{K}{constK:int(-u)e^F}=K(\kappa,\omega_Y,H.p_0)>0$ such that
	\begin{align*}
		\int_Y(-u) e^F\mu_\delta\leq K_{\ref{constK:int(-u)e^F}}.
	\end{align*}
\end{lemma}
\begin{proof}
		Since $u$ is $\kappa\omega_Y$-psh,  there exists a constant $\alpha>0$ depending on $\alpha([\kappa\omega_Y])$ such that
	\begin{align*}
		\int_Ye^{-\alpha u}\omega_Y^n\leq C_Y.
	\end{align*}
	Then by \eqref{eq:Lp-int-mu-delta} and H\"older inequality, we have
	\begin{align*}
		\int_Y e^{-\alpha u}\mu_\delta\leq K.
	\end{align*}
	By \cite[Proposition~2.10]{BBEGZ19}, we obtain
	\begin{align*}
		\int_Y(-\alpha u)\,e^F\mu_\delta
		\leq H_{\eps,\mu_\delta}(u)+\log\int_Ye^{-\alpha u}\mu_\delta
		\leq K.
	\end{align*}
\end{proof}

\begin{proof}[Proof of \textup{Theorem~\ref{thm:unif-C0}}]

	We consider the following auxiliary
	\begin{align*}
		V_\eps^{-1}(\omega_\eps+\ddc\psi)^n
		=\frac{\sqrt{1+F^2}-u}{A_u}e^F\mu_\delta, \quad
		\sup\psi=0,
	\end{align*}
	where
	\begin{align*}
		A_u\coloneqq\int_Y(\sqrt{1+F^2}-u)e^F\mu_\delta.
	\end{align*}
	Note that $A_u\geq 1$ and by $\sqrt{1+x^2}\leq 1+|x|$,
	\begin{align*}
		A_u
		\leq \int_Y(1+|F|)e^F\mu_\delta+K_{\ref{constK:int(-u)e^F}}
		\leq \constant{K}{constK:A_u-bound}.
	\end{align*}
	In particular, 
	\begin{align*}
		\frac{\sqrt{1+F^2}-u}{A_u}
		\geq \frac{1}{K_{\ref{constK:A_u-bound}}}>0.
	\end{align*}
	One chooses a constant $B\geq K_{\ref{constK:theta_eps,t-bound}}+2$. Then for any $\lambda\in(0,1)$,
	\begin{align*}
		\theta_{\eps,t}+(B-\lambda)\omega_\eps
		\geq (B-\lambda-K_{\ref{constK:theta_eps,t-bound}})\omega_\eps\geq0.
	\end{align*}
	We consider the test function $Q\coloneqq F+\lambda\psi-Bu$. Let $z\in Y$ be a maximum point.
	At $z_0$,
	\begin{align*}
		0\geq \Delta_{\omega'}Q
		=&\Delta_{\omega'}F+ \lambda\Delta_{\omega'}\psi
		-B \Delta_{\omega'}u \\
		=&(\tr_{\omega'}\theta_{\eps,t}-c_{\eps,t})+\lambda\tr_{\omega'}(\omega_{\eps,\psi}-\omega_\eps)-B\tr_{\omega'}(\omega'-\omega_\eps) \\
		\geq&-Bn-c_{\eps,t}+\lambda\tr_{\omega'}\omega_{\eps,\psi}.
	\end{align*}
	Note that
	\begin{align*}
		\tr_{\omega'}\omega_{\eps,\psi}
		\geq n\left(\frac{\omega_{\eps,\psi}^n}{\omega'^n} \right)^{1/n}
		=n\left(\frac{\sqrt{1+F^2}-u}{A_u} \right)^{1/n}.
	\end{align*}
	Combining with the above inequality, we have
	\begin{align*}
		(\sqrt{1+F^2}-u)(z)
		\leq  A_u(B+K_{\ref{constK:theta_eps,t-bound}})^n\lambda^{-n}
		\leq K\lambda^{-n},
	\end{align*}
	which implies
	\begin{align*}
		|F(p)|\leq K\lambda^{-n}, \quad -u(p)\leq K\lambda^{-n}.
	\end{align*}
	Thus
	\begin{align*}
		Q\leq Q(p)\leq |F(p)|-Bu(p)\leq K\lambda^{-n}.
	\end{align*}
	It follows that
	\begin{align}
	\label{eq:F-leq-lambda-psi}
		F\leq K\lambda^{-n}-\lambda\psi.
	\end{align}
	\medskip
	
      By \eqref{eq:Lp-int-mu-delta} and H\"older inequality, for some $p>1$, we have
	\begin{align*}
		\int_Y e^{pF}f_\delta^p\,\omega_Y^n
		\leq& K_\lambda\int_Ye^{-p\lambda\psi}f_\delta^{p-1}\,\mu_\delta \\
		\leq& K_\lambda\left(\int_Ye^{-qp\lambda\psi}\mu_\delta \right)^{1/q} \left(\int_Yf_\delta^{q^*(p-1)}\mu_\delta \right)^{1/q^*}
		\leq K_{\lambda},
	\end{align*}
	where we choose $p-1\ll1$ such that $q^*(p-1)<p_0-1$ and $\lambda\leq \frac{\alpha}{pq}$.
	Then by Ko{\l}odziej's $L^\infty$-estimate (\cite{Kol98}), we have
	\begin{align*}
		\|u\|_{L^\infty}\leq K_H.
	\end{align*} 
	
	\noindent{\bf Bounds of $\psi$ and $\sup F$.}
	\medskip
	
	One has 
	\begin{align*}
		\sqrt{1+F^2}-u
		\leq K_H(1+|F|) \quad \text{ and }\quad
		(\sqrt{1+F^2}-u)e^F\leq K_H(1+e^{2F}).
	\end{align*}
	By \eqref{eq:Lp-int-mu-delta} and H\"older inequality, for some $p'>1$, then
	\begin{align}
		\int_Y\left(\frac{\sqrt{1+F^2}-u}{A_u} e^Ff_\delta \right)^{p'}\omega_Y^n
		\leq& K_H(\int_Yf_\delta^{p'}\omega_Y^n+\int_Ye^{2p'F}f_\delta^{p'}\omega_Y^n) \nonumber \\ 
		\leq& K_H+K_H\int_Y e^{2p'F}f_\delta^{p'-1}\mu_\delta \nonumber   \\
		\leq& K_H+K_H\left(\int_Ye^{2q'p'F}\mu_\delta \right)^{1/q'}\left(\int_Yf_\delta^{q'^*(p'-1)}\mu_\delta \right)^{1/q'^*} \nonumber \\
		\leq& K_H+K_H \left(\int_Ye^{2q'p'F}\mu_\delta \right)^{1/q'},
		\label{eq:Lp'-int-density-psi}
	\end{align}
	where we choose $p'-1\ll1$ such that $q'^*(p'-1)+1<p_0$.
	Note that
	\begin{align*}
		\int_Ye^{2q'p'F}\mu_\delta 
		\leq K_\lambda \int_Ye^{-2p'p'\lambda\psi}f_\delta\omega_Y^n
		\leq K_\lambda\left( \int_Ye^{-2p_0^*p'p'\lambda\psi}\omega_Y^n\right)^{1/p_0^*}\|f_\delta\|_{L^{p_0}}
		\leq K_H.
	\end{align*}
	Combining with \eqref{eq:Lp'-int-density-psi}, we obtain
	\begin{align*}
		\int_Y\left(\frac{\sqrt{1+F^2}-u}{A_u} e^Ff_\delta \right)^{p'}\omega_Y^n
		\leq K_H.
	\end{align*}
	Again by Ko{\l}odziej's $L^\infty$-estimate (\cite{Kol98}), we have
	\begin{align*}
		\|\psi\|_{L^\infty}\leq K_H.
	\end{align*} 
	By \eqref{eq:F-leq-lambda-psi}, we obtain
	\begin{align*}
		\sup F\leq K_{\lambda,H}
	\end{align*}
	
	\noindent{\bf Lower bound for $F$.}
	\medskip
	
	We choose another auxiliary equation
	\begin{align*}
		V_\eps^{-1}(\omega_\eps+\ddc\psi')^n=f_\delta\omega_Y^n, \quad
		\sup \psi'=0.
	\end{align*}
	By \eqref{eq:Lp-int-mu-delta} and Ko{\l}odziej's $L^\infty$-estimate, we have $\|\psi\|_{L^\infty}\leq K$.

	We consider the test function $Q'\coloneqq F+Bu-\psi'$. Let $z'\in Y$ be a minimum point of $Q'$.
	By maximum principle, we have,  at $z'$,
	\begin{align*}
		0\leq \Delta_{\omega'}Q'
		=&(\tr_{\omega'}\theta_{\eps,t}-c_{\eps,t})+B\tr_{\omega'}(\omega'-\omega_\eps) 
		-\tr_{\omega'}(\omega_{\eps,\psi'}-\omega_\eps) \\
		\leq&(K_{\ref{constK:theta_eps,t-bound}}-(B-1))\tr_{\omega'}\omega_\eps-c_{\eps,t}+Bn-\tr_{\omega'}\omega_{\eps,\psi'}
		\leq n(K_{\ref{constK:theta_eps,t-bound}}+B)-\tr_{\omega'}\omega_{\eps,\psi'}.
	\end{align*}
	Note that
	\begin{align*}
		\tr_{\omega'}\omega_{\eps,\psi'}
		\geq n\left(\frac{\omega_{\eps,\psi'}^n}{\omega'^n} \right)^{1/n}
		=n\left(\frac{V_\eps f_\delta\omega_Y^n}{V_\eps e^F\mu_\delta} \right)^{1/n}
		=ne^{-F/n}.
	\end{align*}
	Combining with the above inequality, we have
	\begin{align*}
		F(z')\geq-n\log(K_{\ref{constK:theta_eps,t-bound}}+B)\geq-K.
	\end{align*}
	Then
	\begin{align*}
		F+Bu-\psi'\geq F(z')+Bu(z')-v(z')
		\geq-K-B\|u\|_{L^\infty}-\|v\|_{L^\infty},
	\end{align*}
	which implies that
	\begin{align*}
		F\geq-K_H.
	\end{align*}
\end{proof}

\subsection{Uniform $W^{2,p}$-estimate}
In this subsection, we prove the uniform $W^{2,p}$-estimate for the equation \eqref{eq:eps_twist_cscK}.

\bigskip
We fix a smooth log resolution $\pi:Y\rightarrow X$ with exceptional divisor $E \subset Y$  as in \cite[Section~4.1]{GSS26}, that is, $Y$ admits an effective exceptional divisor whose negative is relatively ample, denoted by $E_0$. 
Then $A\coloneqq \pi^*L-\gamma E_0$ is $\Bbb Q$-ample for some small $\gamma>0$ rational. 
Thus there exists a smooth qpsh function $\psi'$ such that $\omega_Y+\gamma c_1(E_0,h_0)=\pi^*\omega+\ddc\psi'$.
	By Poincar\'e-Lelong formula, $\ddc\log|s_{E_0}|^2_{h_0}=[E_0]-c_1(E_0,h_0)$. 
	Thus, for $\gamma$ small enough, 
	\begin{align}
	\label{def_of_psi}
		\omega_Y+\gamma[E_0]=\hat\omega+\ddc(\psi'+\gamma\log|s_{E_0}|_{h_0}^2)
	\end{align}
	 is a big current. So $\psi_0:=\psi'+\gamma\log|s_{E_0}|_{h_0}^2$ is a $\pi^*\omega_X$-qpsh function.
	 Write $E_0=\sum_\ell d_\ell E_\ell$ for some $d_\ell>0$. We may assume
	 \begin{align*}
	 	\psi_0=\psi'+\gamma(\sum_\ell d_\ell\log|s_\ell|_{h_\ell}^2)\leq-1.
	 \end{align*}
After rescaling the $h_i$ by positive constants, we also assume that $|s_i|_{h_i}^2\leq e^{-2}$.
In addition,
\begin{align*}
	|\Delta_{\omega_Y}\psi_0|=|\tr_{\omega_Y}\ddc\psi_0|\leq K, \text{ on }Y\setminus E,\quad \text{ and }
	|\partial\psi_0|_{\omega_Y}\leq K\sum_i|s_i|_{h_i}^{-1}.
\end{align*}
Note that $\psi_0\leq C+ \gamma d_i\log|s_i|_{h_i}^2$. Then there exists $K>0$ such that
\begin{align}
\label{eq:psi_0-f_delta}
	K^{-1}e^{k_0\psi_0}\leq f_\delta\leq Ke^{k_1\psi_0},
\end{align}
where
\begin{align*}
	k_0\coloneqq\sum_{a_i>0}\frac{a_i}{\gamma d_i}\geq0, \qquad
	k_1\coloneqq\sum_{a_i<0}\frac{a_i}{\gamma d_i}\leq0.
\end{align*}

Write $a_i^+\coloneqq \max(a_i,0)$ and $a_i^-\coloneqq\max(-a_i,0)$. We denote $r_i\coloneqq|s_i|^2_{h_i}$ and $\psi_\delta^\pm\coloneqq\sum_i a_i^\pm\log(r_i+\delta)$
\begin{lemma} As the above notations, there exists  $K>0$ such that
\label{lem:lowerbound-ddclog-f_delta}
	\begin{align*}
		\Delta_{\omega_Y}\log(f_\delta)\geq -K-\Delta_{\omega_Y}\psi_\delta^- .
	\end{align*}
\end{lemma}
\begin{proof}
Recall that $\log f_\delta= \log h_0+\psi_\delta^+-\psi_\delta^-$.
Note that
	\begin{align*}
		\ddc\log(r_i+\delta)
		=\frac{r_i}{r_i+\delta}\ddc\log r_i
		+\frac{\delta}{r_i(r_i+\delta)^2}\partial r_i\wedge\bp r_i
	\end{align*}
	and
	\begin{align*}
		\left|\frac{r_i}{r_i+\delta}\ddc\log r_i\right|
		\leq \frac{r_i}{r_i+\delta}|\Theta_i|
		\leq K\omega_Y.
	\end{align*}
	Then one has
	\begin{align*}
		\ddc\log(r_i+\delta)
		\geq -K\omega_Y,
	\end{align*}
	which yields that $\ddc\psi_\delta^\pm\geq-\constant{K}{constK;qpsh-psi_delta^pm}\omega_Y$ for some $K_{\ref{constK;qpsh-psi_delta^pm}}>0$. Thus we obtain
	\begin{align*}
		\Delta_{\omega_Y}\log(f_\delta)\geq -K-\Delta_{\omega_Y}\psi_\delta^- .
	\end{align*}
\end{proof}

\begin{theorem}
	\label{thm:unif-W2p}
	Suppose $(u, F)$ is a solution of \eqref{eq:eps_twist_cscK}. For any $p\geq 1$, there exist uniform constants $a,b>1$ and $K_p>0$ such that
	\begin{align*}
		\int_Y e^{\big(ab(2p+1)+k_0-\frac{nc_1}{n-1}\big)\psi_0}(\tr_{\omega_Y}\omega')^{2p+1+\frac{1}{n-1}}\omega_Y^n
		\leq K_p=K(p,\|u\|_{L^\infty}, \|F\|_{L^\infty},\omega_Y).
	\end{align*}
\end{theorem}

\begin{lemma}
\label{lem:tr-relation}
	There exists a constant $K>0$ such that
	\begin{align*}
		\tr_{\omega'}\omega_Y\geq Ke^{-\frac{F}{n-1}}e^{-\frac{k_1}{n-1}\psi_0} (\tr_{\omega_Y}\omega')^{\frac{1}{n-1}}.
	\end{align*}
\end{lemma}
\begin{proof}
	Note that, by \eqref{eq:eps_twist_cscK},
	\begin{align*}
		\tr_{\omega_Y}\omega'
    \leq \frac{1}{(n-1)!}(\tr_{\omega'}\omega_Y)^{n-1}\frac{\omega'^n}{\omega_Y^n}
    =\frac{1}{(n-1)!}(\tr_{\omega'}\omega_Y)^{n-1}e^FV_\eps f_\delta.
	\end{align*}
	Then we have
	\begin{align*}
		\tr_{\omega'}\omega_Y\geq Ke^{-\frac{F}{n-1}}f_\delta^{-\frac{1}{n-1}} (\tr_{\omega_Y}\omega')^{\frac{1}{n-1}}.
	\end{align*}
	Combining with \eqref{eq:psi_0-f_delta}, we finish the proof.
\end{proof}

\begin{proof}[Proof of \textup{Theorem~\ref{thm:unif-W2p}}]
	Since $\psi_0\in\PSH_{\kappa\omega_Y}$, then there exists a family of $\psi_j\in\Psh_{\kappa\omega_Y}\cap C^\infty$ which decreasingly converges to $\psi_0$.
	
	We consider test functions
	\begin{align*}
		w\coloneqq e^{ab\psi_0-a(F+bu)+\psi_\delta^-}\tr_{\omega_Y}\omega^\prime \quad
		\text{ and }\quad
		w_j\coloneqq e^{ab\psi_j-a(F+bu)+\psi_\delta^-}\tr_{\omega_Y}\omega^\prime,
	\end{align*}
	where $a,b>1$ are uniform constants determined later.
	Then $w_j$ is a family of smooth qpsh function, decreasingly converging to the bounded qpsh function $w$.

    One has
	\begin{align}
		\Delta_{\omega^\prime}w_j
		=\Delta_{\omega^\prime} e^{\log w_j}
		=&e^{\log w_j}\Delta_{\omega^\prime}\log w_j
		 +w_j\nabla^{\omega'}(\log w_j)\cdot_{\omega'}\nabla^{\omega'}(\log w_j) \nonumber \\
		=&w_j\Delta_{\omega'}(ab\psi_j-a F-abu+\psi_\delta^-) +w_j\Delta_{\omega^\prime}\log\tr_{\omega_Y}\omega^\prime +w_j^{-1}|\nabla^{\omega'} w_j|_{\omega^\prime}^2 .
		\label{eq:Lap_phi e logu}
	\end{align}
	By Yau's Schwarz lemma \cite{Yau78b}  (see also e.g. \cite[Chap 3, Lemma~3.7]{Sze14}),
\begin{align}
\label{eq:Yau-inequ}
	\Delta_{\omega'} \log\tr_{\omega_Y}\omega'
	\geq -B\tr_{\omega'}\omega_Y-\frac{\tr_{\omega_Y}\Ric(\omega')}{\tr_{\omega_Y}(\omega')},
\end{align}
where $-B$ is a lower bound for the bisectional curvature of $\omega_Y$. 
By \eqref{eq:eps_twist_cscK} and Lemma~\ref{lem:lowerbound-ddclog-f_delta}, one has
\begin{align*}
	\tr_{\omega_Y}\Ric(\omega')
	=\tr_{\omega_Y}(-\ddc F-\ddc\log f_\delta+\Ric(\omega_Y))
	\leq -\Delta_{\omega_Y}F+\Delta_{\omega_Y}\psi_\delta^-+\constant{K}{constK:bound-f_delta+Ric(omega)},
\end{align*}
where $K_{\ref{constK:bound-f_delta+Ric(omega)}}>0$ is a constant such that $\Ric(\omega_Y)\leq K_{\ref{constK:bound-f_delta+Ric(omega)}}\omega_Y$.
Since $(\tr_{\omega_Y}\omega')^{-1}\leq \frac{1}{n^2}\tr_{\omega'}\omega_Y$, then we have
\begin{align}
\label{eq:Yau-inequ-simplify}
	\Delta_{\omega'} \log\tr_{\omega_Y}\omega'
	\geq -(B+K_{\ref{constK:bound-f_delta+Ric(omega)}})\tr_{\omega'}\omega_Y
	+\frac{\Delta_{\omega_Y}F-\Delta_{\omega_Y}\psi_\delta^-}{\tr_{\omega_Y}(\omega')}.
\end{align}
Since $K_{\ref{constK;qpsh-psi_delta^pm}}\omega_Y+\ddc\psi_\delta^-\geq0$ for , then 
\begin{align*}
	0\leq K_{\ref{constK;qpsh-psi_delta^pm}}\omega_Y+\ddc\psi_\delta^-
	\leq \tr_{\omega'}(K_{\ref{constK;qpsh-psi_delta^pm}}\omega_Y+\ddc\psi_\delta^-)\omega'.
\end{align*}
Applying $\tr_{\omega_Y}$ on both sides,
\begin{align*}
	0\leq \tr_{\omega_Y}(K_{\ref{constK;qpsh-psi_delta^pm}}\omega_Y+\ddc\psi_\delta^-)
	\leq \tr_{\omega'}(K_{\ref{constK;qpsh-psi_delta^pm}}\omega_Y+\ddc\psi_\delta^-)\tr_{\omega_Y}\omega'.
\end{align*}
and so
\begin{align}
\label{eq:Lap_psi-}
	\Delta_{\omega_Y}\psi_\delta^-
	\leq \tr_{\omega'}(K_{\ref{constK;qpsh-psi_delta^pm}}\omega_Y+\ddc\psi_\delta^-)\tr_{\omega_Y}\omega'.
\end{align}
Combining with \eqref{eq:Yau-inequ-simplify}, we have
\begin{align}
\label{eq:Yau-inequ-simplify-final}
	\Delta_{\omega'} \log\tr_{\omega_Y}\omega'
	\geq -(B+K_{\ref{constK:bound-f_delta+Ric(omega)}}+K_{\ref{constK;qpsh-psi_delta^pm}})\tr_{\omega'}\omega_Y
	+\frac{\Delta_{\omega_Y}F}{\tr_{\omega_Y}(\omega')}
	-\Delta_{\omega'}\psi_\delta^-.
\end{align}
	
   On the other hand, by \eqref{eq:eps_twist_cscK}, we have
	\begin{align}
		\Delta_{\omega^\prime}(F+bu)
		=&\tr_{\omega^\prime}\theta_{\eps,t}-c_{\eps,t}
		  +b\tr_{\omega^\prime}(\omega^\prime-\omega_\epsilon) \nonumber \\
		\leq&(bn-c_{\eps,t}) + \kappa K_{\ref{constK:theta_eps,t-bound}}\tr_{\omega^\prime}\omega_Y-b\tr_{\omega'}\omega_\epsilon. \label{Lap_phi(tilde F+bphi)}
	\end{align}
	Note that
	\begin{align*}
		\omega_\epsilon\wedge\omega'^{n-1}
		=(\hat\omega+\epsilon\omega_Y)\wedge\omega'^{n-1}
		=&(1+\epsilon)\omega_Y\wedge\omega'^{n-1}+(\hat\omega-\omega_Y)\wedge\omega'^{n-1} \\
		=&(1+\epsilon)\omega_Y\wedge\omega'^{n-1}+(\gamma[E_0]-\ddc\psi_0)\wedge\omega'^{n-1} \\
		\geq&(1+\epsilon)\omega_Y\wedge\omega'^{n-1} -\ddc\psi_0\wedge\omega'^{n-1}.
	\end{align*}
	Together with \eqref{eq:Lap_phi e logu} and \eqref{eq:Yau-inequ-simplify-final},  we have
	\begin{align}
		&n\ddc w_j\wedge\omega'^{n-1} \\
		\geq& e^{ab\psi_j-a(F+bu)+\psi_\delta^-+\log\tr_{\omega_Y}\omega'}n\ddc(ab\psi_j-a( F+bu)+\log\tr_{\omega_Y}\omega')\wedge\omega'^{n-1} \nonumber \\
		\geq& e^{ab\psi_j-a(F+bu)+\psi_\delta^-}
		\Big(a(c_{\eps,t}-bn)\tr_\omega\omega^\prime
		+\Big((1+\epsilon)ab-a\kappa K_{\ref{constK:theta_eps,t-bound}}-B-K_{\ref{constK:bound-f_delta+Ric(omega)}}\Big)
		\tr_{\omega_Y}\omega^\prime\tr_{\omega^\prime}\omega_Y
		+\Delta_{\omega_Y} F \Big)\omega'^n \nonumber \\ 
		&+w_jnab(\ddc\psi_j-\ddc\psi_0)\wedge\omega'^{n-1} \nonumber \\
		\geq&aw_j(c_{\eps,t}-bn)\omega'^n +\frac{ab}{2}Kw_je^{-\frac{F}{n-1}}e^{-\frac{k_1}{n-1}\psi_0}(\tr_\omega\omega')^{\frac{1}{n-1}}\cdot\omega'^n +e^{ab\psi_j-a( F+bu)+\psi_\delta^-}\Delta_\omega F\cdot\omega'^n \nonumber \\
		&+w_jnab(\ddc\psi_j-\ddc\psi_0)\wedge\omega'^{n-1}.
		\label{eq:tr of ddc w_j geq ...}
	\end{align}
	where we choose $b\gg1$ such that $ab-a\kappa K_{\ref{constK:theta_eps,t-bound}}-B-K_{\ref{constK:bound-f_delta+Ric(omega)}}\geq\frac{ab}{2}$ and we have used Lemma~\ref{lem:tr-relation}.
	
	Note that 
	\begin{align*}
		\frac{1}{2p+1}\Delta_{\omega^\prime} w_j^{2p+1}
		=&w_j^{2p}\Delta_{\omega^\prime} w_j +2pw_j^{2p-1}|\nabla^{\omega^\prime} w_j|_{\omega^\prime}^2 \\
		\geq& w_j^{2p}\Delta_{\omega^\prime}w_j
		   +2pw_j^{2p-2}e^{ab\psi_j-a(F+bu)+\psi_\delta^-}|\nabla^{\omega_Y} w_j|_{\omega_Y}^2,
	\end{align*}
	where we have used
	\begin{align*}
		|\nabla^{\omega^\prime} w_j|_{\omega^\prime}^2\tr_{\omega_Y}\omega^\prime\geq|\nabla^{\omega_Y} w_j|_{\omega_Y}^2.
	\end{align*}
	 Together with \eqref{eq:tr of ddc w_j geq ...}, we have
	\begin{align}
		0=&\frac{1}{2p+1}\int_Y\Delta_{\omega^\prime} w_j^{2p+1}\omega^{\prime n} \nonumber \\
		\geq& 2p\int_Yw_j^{2p-2}|\nabla^{\omega_Y} w_j|_{\omega_Y}^2e^{ab\psi_j-a(F+bu)+\psi_\delta^-+ F+\log(V_\eps f_\delta)}\omega_Y^n 
		  +a(c_{\eps,t}-bn)\int_Yw_j^{2p+1}\omega'^n  \nonumber \\
		  &+\int_Yw_j^{2p}e^{ab\psi_j-a( F+bu)+\psi_\delta^-+F+\log(V_\eps f_\delta)}\Delta_{\omega_Y} F\omega_Y^n  
		  +\frac{ab}{2}K\int_Yw_j^{2p+1}(\tr_{\omega_Y}\omega^{\prime})^{\frac{1}{n-1}} e^{-\frac{1}{n-1} F}e^{-\frac{k_1}{n-1}\psi_0}\omega'^n \nonumber\\
		 &+nab\int_Yw_j^{2p+1}\ddc(\psi_j-\psi_0)\wedge\omega'^{n-1}.
		   \label{eq:integ of Lap_u w_j^{2p+1}}
	\end{align}
	By Lemma~\ref{lem:int-I-term} and choosing $a$ big enough (eg. $a=kp+1$ for $k\gg 1$) such that $2p\geq\frac{p^2}{a-1}$, one has
	\begin{align}
		0\geq&nab\int_Yw_j^{2p+1}\ddc(\psi_j-\psi)\wedge\omega'^{n-1} 
		   +\frac{nab}{a-1}\int_Yw_j^{2p}e^{G_j}(\ddc\psi_j-\ddc\psi)\wedge \omega^{n-1} \nonumber \\
		 &-K\int_Yw_j^{2p+1}e^{ F}\omega'^n 
		 -\frac{K}{a-1}\int_Yw_j^{2p}e^{G_j}\omega_Y^n
		  +\frac{ab}{2}K\int_Yw_j^{2p+1}(\tr_{\omega_Y}\omega^{\prime})^{\frac{1}{n-1}} e^{-\frac{1}{n-1} F}e^{-\frac{k_1}{n-1}\psi_0}\omega'^n,
		  \label{eq:integration of w_j}
	\end{align}
	where $K>0$ is a uniform constant depending on $a, b, c_{\eps,t}, n$.	
	
	Note that for any $\epsilon>0$,
	 $w_j\searrow w$ bounded, by Hartog's Lemma, $w_j\rightarrow w$ in $L^q$ for any $q\geq1$.
	 Similarly, $e^{G_j}\rightarrow e^G$ in $L^q$ for any $q\geq1$, where $G\coloneqq ab\psi_0+(1-a) F-abu+\psi_\delta^-+\log(V_\eps f_\delta)$. Then
	\begin{align*}
		\int_Yw_j^{2p+1}e^{F}\omega'^n\rightarrow
		\int_Yw^{2p+1}e^{ F}\omega'^n,\quad
		\int_Yw_j^{2p}e^{G_j}\omega_Y^n
		\rightarrow
		\int_Yw^{2p}e^{G}\omega_Y^n
	\end{align*}
	and
	\begin{align*}
		\int_Yw_j^{2p+1}(\tr_{\omega_Y}\omega^{\prime})^{\frac{1}{n-1}} e^{-\frac{1}{n-1} F}e^{-\frac{k_1}{n-1}\psi_0}\omega'^n
		\rightarrow
		\int_Yw^{2p+1}(\tr_{\omega_Y}\omega^{\prime})^{\frac{1}{n-1}} e^{-\frac{1}{n-1} F}e^{-\frac{k_1}{n-1}\psi_0}\omega'^n.
	\end{align*}
	Note that
	\begin{align*}
		w_j^{2p+1}=e^{(2p+1)ab\psi_j-(2p+1)a( F+bu)+(2p+1)\psi_\delta^-+(2p+1)\log\tr_\omega\omega'}
	\end{align*}
	and for any $\epsilon>0$, $(2p+1)ab\psi_j-(2p+1)a( F+bu)+(2p+1)\psi_\delta^-+(2p+1)\log\tr_\omega\omega'$ is $\kappa_\epsilon\omega$-psh function for some $\kappa_\epsilon>0$.
	Thus, for any $\epsilon>0$, it follows from Lemma \ref{lem:convergence lemma} that 
	\begin{align*}
		\int_Yw_j^{2p+1}\ddc(\psi_j-\psi)\wedge\omega'^{n-1}\rightarrow0
		\quad\text{ and }\quad\int_Yw_j^{2p}e^{G_j}(\ddc\psi_j-\ddc\psi)\wedge \omega^{n-1}\quad\text{as }j\rightarrow\infty.
	\end{align*}
	Combining the above convergences and \eqref{eq:integration of w_j}, we have
	\begin{align} 
	\label{eq:tr-L2p+1+1/n-1}
		  \frac{ab}{2}K\int_Yw^{2p+1}(\tr_{\omega_Y}\omega^{\prime})^{\frac{1}{n-1}} e^{-\frac{1}{n-1} F}e^{-\frac{k_1}{n-1}\psi_0}\omega'^n		  \leq K\int_Yw^{2p+1}e^{ F}\omega'^n 
		 +K\int_Yw^{2p}e^{G}\omega_Y^n.
	\end{align}
	We denote 
	\begin{align*}
		\nu\coloneqq e^{\big(1-a(2p+1)-\frac{1}{n-1}\big)F +\big(ab(2p+1)-\frac{k_1}{n-1}\big)\psi_0-(2p+1)abu+(2p+1)\psi_\delta^-}V_\eps f_\delta\omega_Y^n
	\end{align*} 
	Then, since $\psi_\delta^-$ is bounded from above,
	\begin{align*}
		\nu(Y)
		\leq K\int_Y e^{\big(ab(2p+1)-\frac{k_1}{n-1}\big)\psi_0}f_\delta\omega_Y^n
		\leq K\int_Y e^{\big(ab(2p+1)+\frac{n-2}{n-1}k_1\big)\psi_0}f_\delta\omega_Y^n
		\leq K,
	\end{align*}
	where we have used $ab(2p+1)+\frac{n-2}{n-1}k_1\geq0$ and $\psi_0\leq-1$.
	Then the inequality \eqref{eq:tr-L2p+1+1/n-1} can be re-written as
	\begin{align}
	\label{eq:N-def-ineq}
		N\coloneqq\int_Y(\tr_{\omega_Y}\omega^{\prime})^{2p+1+\frac{1}{n-1}}\,d\nu
		\leq K\int_Y(\tr_{\omega_Y}\omega^{\prime})^{2p+1}e^{\frac{k_1}{n-1}\psi_0}\,d\nu
		+K\int_Y(\tr_{\omega_Y}\omega^{\prime})^{2p}e^{\frac{k_1}{n-1}\psi_0}\,d\nu.
	\end{align} 
	For $k=2p$ or $k=2p+1$, by H\"older inequality with the index $r_k=\frac{2p+1+\frac{1}{n-1}}{k}$, one has
	\begin{align*}
		\int_Y(\tr_{\omega_Y}\omega^{\prime})^{k}e^{\frac{k_1}{n-1}\psi_0}\,d\nu
		\leq \left(\int_Y\tr_{\omega_Y}\omega^{\prime})^{2p+1+\frac{1}{n-1}}d\nu\right)^{\frac{k}{2p+1+\frac{1}{n-1}}}\left(\int_Y e^{r_k^*\frac{k_1}{n-1}\psi_0}d\nu\right)^{1/r_k^*}
		\leq K_pN^{\frac{k}{2p+1+\frac{1}{n-1}}},
	\end{align*}
	where we have used the fact
	\begin{align*}
		\int_Y e^{r_k^*\frac{k_1}{n-1}\psi_0}d\nu
		\leq K\int_Y e^{\big(ab(2p+1)+k_1\frac{r_k^*+n-2}{n-1}\big)\psi_0}\omega_Y^n
		\leq  K_p,
	\end{align*}
	which is due to $\psi_0\leq-1$ and the choice of $b$, such that
	\begin{align*}
		ab(2p+1)+k_1\frac{r_k^*+n-2}{n-1}\geq0.
	\end{align*}
	Hence, combining with \eqref{eq:N-def-ineq} and Young's inequality ($x^\beta\leq \eps x+C_\beta$ for $\beta<1$), we obtain
	\begin{align*}
		N\leq K_p\left(N^{\frac{2p+1}{2p+1+\frac{1}{n-1}}} 
		+N^{\frac{2p}{2p+1+\frac{1}{n-1}}} \right)
		\leq \frac{1}{2}N+K_p,
	\end{align*}
	that is,
	\begin{align*}
		\int_Y(\tr_{\omega_Y}\omega^{\prime})^{2p+1+\frac{1}{n-1}}\,d\nu
		\leq K_p.
	\end{align*}
	Since $\log(r_i+\delta)\geq \log r_i\geq (\gamma d_i)\psi_0-C$, then $\psi_\delta^-\geq(-k_1)\psi_0-C$
	According to the definition of $\nu$, we obtain
		\begin{align*}
			\int_Y(\tr_{\omega_Y}\omega^{\prime})^{2p+1+\frac{1}{n-1}}e^{\big(ab(2p+1)+k_0-\frac{nk_1}{n-1}\big)\psi_0}\omega_Y^n
		\leq K_p.
		\end{align*}
\end{proof}

\begin{lemma}
\label{lem:int-I-term}
We denote $G_j\coloneqq ab\psi_j+(1-a) F-abu+\psi_\delta^-+\log(V_\eps f_\delta)$. Then
	\begin{align*}
		\int_Yw_j^{2p}e^{G_j}\Delta F\omega_Y^n
		\geq& -\frac{p^2}{a-1}\int_Yw_j^{2p-2}|\nabla w_j|^2e^{G_j}\omega_Y^n
		-\frac{ab}{a-1}\int_Yw_j^{2p+1}e^{ F}\omega'^n
		-\frac{K}{a-1}\int_Yw_j^{2p}e^{G_j}\omega_Y^n\\
		&+\frac{nab}{a-1}\int_Yw_j^{2p}e^{G_j}(\ddc\psi_j-\ddc\psi_0)\wedge \omega_Y^{n-1}.
	\end{align*}
\end{lemma}

\begin{proof}
	For simplicity, we denote $\Delta=\Delta_{\omega_Y}$. 
	Note that $G_j\coloneqq ab\psi_j+(1-a) F-abu+\psi_\delta^++\log h_0$.
	Thus 
	\begin{align*}
		F=\frac{1}{1-a}G_j+\frac{ab}{1-a}u-\frac{ab}{1-a}\psi_j-\frac{1}{1-a}(\psi_\delta^++\log h_0).
	\end{align*}
	One has
	\begin{align*}
		I&\coloneqq \int_Yw_j^{2p}e^{G_j}\Delta F\omega_Y^n
		 \\
	     &\ =\frac{1}{1-a}\int_Yw_j^{2p}\Delta G_je^{G_j}\omega_Y^n
		 +\frac{ab}{1-a}\int_Yw_j^{2p}\Delta (u-\psi_j) e^{G_j}\omega_Y^n  
		 -\frac{1}{1-a}\int_Yw_j^{2p}\Delta(\psi_\delta^++\log h_0) e^{G_j}\omega_Y^n \\
		 &\ \eqqcolon I_1+I_2+I_3.
	\end{align*}
	By Stokes' theorem and Young inequality, we obtain
	\begin{align*}
		I_1=&\frac{1}{a-1}\int_Yw_j^{2p}|\nabla G_j|^2e^{G_j}\omega_Y^n
		+\frac{2p}{a-1}\int_Yw_j^{2p-1}\nabla w_j\cdot\nabla G_j e^{G_j}\omega_Y^{n} \\
		\geq&\frac{1}{a-1}\int_Yw_j^{2p}|\nabla G_j|^2e^{G_j}\omega_Y^n
		-\frac{p^2}{a-1}\int_Yw_j^{2p-2}|\nabla w_j|^2e^{G_j}\omega_Y^n 
		-\frac{1}{a-1}\int_Yw_j^{2p}|\nabla G_j|^2e^{G_j}\omega_Y^n \\
		\geq&-\frac{p^2}{a-1}\int_Yw_j^{2p-2}|\nabla w_j|^2e^{G_j}\omega_Y^n.
	\end{align*}
	In addition, since $e^{G_j}(\tr_{\omega_Y}\omega')\omega_Y^n=w_je^F\omega'^n$,
	\begin{align*}
		I_2=&\frac{ab}{1-a}\int_Yw_j^{2p}\tr_{\omega_Y}(\omega^\prime-\omega_\epsilon-\ddc\psi_j) e^{G_j}\omega_Y^n \\
		=&-\frac{ab}{a-1}\int_Yw_j^{2p+1}e^{ F}\omega'^n
		+\frac{ab}{a-1}\int_Yw_j^{2p}\tr_{\omega_Y}(\omega_\epsilon+\ddc\psi_j) e^{G_j}\omega_Y^n \\
		=&-\frac{ab}{a-1}\int_Yw_j^{2p+1}e^{F}\omega'^n
		+\frac{nab}{a-1}\int_Yw_j^{2p}e^{G_j}(\omega_\epsilon+\ddc\psi_0+\ddc\psi_j-\ddc\psi_0)\wedge \omega_Y^{n-1} \\
		\geq&-\frac{ab}{a-1}\int_Yw_j^{2p+1}e^{ F}\omega'^n
		+\frac{nab}{a-1}\int_Yw_j^{2p}e^{G_j}(\ddc\psi_j-\ddc\psi_0)\wedge \omega_Y^{n-1},
	\end{align*}
	where we have used that $\omega_\epsilon+\ddc\psi_0=\hat\omega_X+\epsilon\omega+\ddc\psi_0\geq0$.

	Since $\ddc\psi_\delta^\pm\geq-K_{\ref{constK;qpsh-psi_delta^pm}}\omega_Y$ we have
	\begin{align*}
		I_3=-\frac{1}{1-a}\int_Yw_j^{2p}\Delta(\psi_\delta^++\log h_0) e^{G_j}\omega_Y^n
		\geq -\frac{K}{a-1}\int_Yw_j^{2p}e^{G_j}\omega_Y^n.
	\end{align*}
	Combining those three terms, we finish the proof.
\end{proof}

\begin{lemma}
\label{lem:convergence lemma}
	Let $(Y,\omega_Y)$ be a compact K\"ahler manifold and $\psi_0\in\PSH_{\omega_Y}$. Suppose $\psi_j\in\PSH_{\omega_Y}$ is a decreasing sequence of smooth $\omega_Y$-psh functions converging to $\psi$. Then
	\begin{align*}
		\int_Y e^{\psi_j}\ddc(\psi_j-\psi)\wedge\omega_Y^{n-1}\rightarrow0.
	\end{align*} 
\end{lemma}
\begin{proof}
	By Hartog's Lemma, $\psi_j\rightarrow\psi$ in $L^p$ for $p\geq1$.
	Thus, one has $\ddc\psi_j\rightarrow\ddc\psi$ in the sense of current.
	The following formula holds in the sense of current
	\begin{align*}
		\ddc e^\psi
		=&e^\psi\ddc\psi+e^\psi\d\psi\wedge\d^c\psi \\
		=&e^\psi\ddc\psi+4\d e^{\frac{\psi}{2}}\wedge \d^ce^{\frac{\psi}{2}}.
	\end{align*}
	Note that
	\begin{align*}
		&\int_Ye^{\psi_j}\ddc\psi_j\wedge\omega_Y^{n-1}-\int_Ye^{\psi}\ddc\psi\wedge\omega_Y^{n-1} \\
		=&\int_Y\ddc e^{\psi_j}\wedge\omega_Y^{n-1}-\int_Y\ddc e^{\psi}\wedge\omega_Y^{n-1} 
		  +4\int_Y\d e^{\frac{\psi}{2}}\wedge\d^ce^{\frac{\psi}{2}}\wedge\omega_Y^{n-1}
		  -4\int_Y\d e^{\frac{\psi_j}{2}}\wedge\d^ce^{\frac{\psi_j}{2}}\wedge\omega_Y^{n-1} \\
		  \rightarrow&0,
	\end{align*}
	since $e^{\psi_j}$ decreasingly converges to a bound qpsh function $e^{\psi}$.
	
	In addition, by dominated convergence theorem,
	\begin{align*}
		&\int_Y(e^{\psi}-e^{\psi_j})\ddc\psi\wedge\omega_Y^{n-1} \\
		=&\int_Y(1-e^{\psi_j-\psi})e^\psi\ddc\psi\wedge\omega_Y^{n-1} \\
		=&\int_Y(1-e^{\psi_j-\psi})\ddc e^\psi\wedge\omega_Y^{n-1} 
		-4\int_Y(1-e^{\psi_j-\psi})\d e^{\frac{\psi}{2}}\wedge\d^ce^{\frac{\psi}{2}}\wedge\omega_Y^{n-1}  \\
		\rightarrow&0.
	\end{align*}
	Then we have
	\begin{align*}
		&\int_X e^{\psi_j}\ddc(\psi_j-\psi)\wedge\omega_Y^{n-1} \\
		=&\int_Ye^{\psi_j}\ddc\psi_j\wedge\omega_Y^{n-1}-\int_Ye^{\psi}\ddc\psi\wedge\omega_Y^{n-1}
		 +\int_Ye^{\psi}\ddc\psi\wedge\omega_Y^{n-1}-\int_Ye^{\psi_j}\ddc\psi\wedge\omega_Y^{n-1} \\
		 \rightarrow&0.
	\end{align*}
\end{proof}

By the well-known relation between $\tr_{\omega_Y}\omega'$ and $\tr_{\omega'}\omega_Y$, an immediate application of Theorem~\ref{thm:unif-W2p} is
\begin{corollary}
	For any $p>0=1$, there exist uniform constants $p_0>0$ and $K_p>0$ such that
	\begin{align*}
		\int_Y e^{p_0\psi_0}(1+\tr_{\omega_Y}\omega'+\tr_{\omega'}\omega_Y)^p\omega_Y^n\leq K_p, \quad
		\int_Y e^{p_0\psi_0}(1+\tr_{\omega_Y}\omega'+\tr_{\omega'}\omega_Y)^p\omega'^n\leq K_p.
	\end{align*}
\end{corollary}

\subsection{Uniform local higher order estimate}
In this subsection, we prove the uniform local higher order estimate for the equation \eqref{eq:eps_twist_cscK}.
We follow Chen-Cheng's local higher order estimate (\cite[Section~6]{CC17}).
\bigskip

Let $B_1\Subset Y\setminus E$ be a coordinate ball.
The equation \eqref{eq:eps_twist_cscK} can be re-written as
\begin{align*}
	\log\det(\omega')=F+\log(f_\delta\det(\omega_Y)), \quad
	\Delta_{\omega'}F=\tr_{\omega'}\theta_{\eps,t}-c_{\eps,t}.
\end{align*}
Choose local potentials $\ell_\eps,\ell$ on $B_1$ such that $\ddc\ell_\eps=\omega_\eps$ and $\ddc\ell=\theta_{\eps,t}$.
Then 
\begin{align}
\label{eq:bound-local-potential}
	\|\ell_\eps\|_{C^m(B_1,\omega_Y)}+\|\ell\|_{C^m(B_1,\omega_Y)}\leq K.
\end{align}
 We denote $F'\coloneqq F-\ell$ and $u'\coloneqq u+\ell_\eps$.
Then on $B_1$, the equation can be simplified as
\begin{align}
\label{eq:local-eq-simp}
	\log\det(u'_{i\bar j})=F'+f', \qquad \Delta_{u'}F'=-c_{\eps,t},
\end{align} 
where $f'\coloneqq \log(f_\delta\det(\omega_Y))+\ell$ with $\|f'\|_{C^m(B_1,\omega_Y)}\leq K$. 

\begin{proposition}
\label{prop:CC-local-estimate}
	As above notations, there exists a uniform constant $K>0$ such that
	\begin{align}
	\label{eq:local-CC-C2}
		|\nabla F|_{\omega_Y}\leq K, \qquad K^{-1}\omega_Y\leq\omega'\leq K\omega_Y \quad\text{on } B_{\frac{1}{2}}.
	\end{align}
	Moreover, for any $m\geq 2$, there exists a uniform constant $K_m>0$ such that
	\begin{align*}
		\|u\|_{C^m(B_{1/4},\omega_Y)}+\|F\|_{C^m(B_{1/4},\omega_Y)}\leq K_m.
	\end{align*}
\end{proposition}
\begin{proof}
	By Chen-Cheng's local $C^2$-estimate (\cite[Proposition~6.1]{CC17}), on $B_{1/2}$, we have
	\begin{align*}
		|\nabla F'|_{\omega_Y}\leq K, \qquad K^{-1}\leq u'_{i\bar j}\leq K. 
	\end{align*}
	Thus the equation \eqref{eq:local-eq-simp} is uniformly elliptic.
	B y Evans-Krylov estimate, then $u'\in C^{2,\gamma}(B_{7/16})$. Applying Schauder's estimate for the second equation, then $F'\in C^{2,\gamma}(B_{3/8})$.
	
	Let $N=2m-2$ and $r_j=\frac{1}{2}-\frac{j}{4N}$ for $0\leq j\leq N$. 
	Applying the bootstrapping and Schauder interior estimate at each ball $B_{r_j}$, then we obtain
	\begin{align*}
		\|u'\|_{C^m(B_{1/4},\omega_Y)}+\|F'\|_{C^m(B_{1/4},\omega_Y)}\leq K_m.
	\end{align*} 
	By \eqref{eq:bound-local-potential}, we finish the proof.
	\end{proof}
	
	\begin{theorem}
	\label{thm:unif-higher-order}
		Suppose $(u, F)$ is a solution of \eqref{eq:eps_twist_cscK} satisfying 
		\begin{align*}
			\|u\|_{L^\infty}+\|F\|_{L^\infty}\leq K_H
		\end{align*}
		Then
		\begin{enumerate}[$(\rm i)$]
			\item For any compact $F\subseteq Y\setminus E$ and any $m>0$, there exists $K_{F,m}>0$ such that 
			\begin{align*}
		  \|u\|_{C^m(F,\omega_Y)}+\|F\|_{C^m(F,\omega_Y)}\leq K_{F,m}.
	\end{align*}
			\item For fixed $\eps,\delta>0$, for any $m>0$, there exists $K_{\eps,\delta,m}>0$ (uniformly for $t\in[0,1]$) such that 
			\begin{align*}
				\|u\|_{C^m(Y,\omega_Y)}+\|F\|_{C^m(Y,\omega_Y)}\leq K_{\eps,H,m}.
			\end{align*}
		\end{enumerate}
	\end{theorem}
	\begin{proof}
		The part (i) is due to Proposition~\ref{prop:CC-local-estimate}.
		For (ii), we can take $\psi_0=0$ in $W^{2,p}$-estimate. Then we obtain the global $W^{2,p}$-estimate.
		Applying Proposition~\ref{prop:CC-local-estimate} on the finitely many balls of a fixed atlas covering $Y$ yields the global $C^m$-estimate, with all constants uniform in $t\in  [0, 1]$ because $F$, $\theta_{\eps,t}$, $c_{\eps,t}$ and all inputs are $t$-uniformly bounded.
	\end{proof}

\section{Proof of the main theorem}
\label{sec:Proof of the main theorem}
In this section, we prove the main theorem. First, we use the Schauder fix point theorem to solve the perturbed equation \eqref{eq:eps_twist_cscK} for $t<1$.
Then according to the uniform estimates in Section~\ref{sec:Uniform esitmates for perturbed equations}, we construct the singular cscK metric.

\subsection{Mabuchi functional of the limit solution}
In this subsection, we prove \ that the limit solution of \eqref{eq:eps_twist_cscK} is the minimizer of Mabuchi functional on $\cH_{\omega}$.
\medskip

Suppose $\eps_k,\delta_k\searrow0$ and $t_k\rightarrow t\in [0,1]$.
For simplicity, we use the conventions  
\begin{align*}
       \omega_k=&\omega_{\eps_k},
       \quad u_k=u_{\eps_k,\delta_k,t_k},  \quad F_k=F_{\eps_k,\delta_k,t_k},
       \quad \omega_k'=\omega_k+\ddc u_k,\quad \mu_k= \mu_{\delta_k},\quad c_k=c_{\eps_k,t_k}, \
         \\
        E_k=&E_{\omega_{\eps_k}}(), \ E_{k,\chi}=E_{\omega_{\eps_k},\chi}(), \quad
         J_{k,\chi}=E_{k,\chi}-\frac{n[\lambda][\omega_k]^{n-1}}{V_k}E_k.
\end{align*}
where $\chi$ is a smooth closed $(1,1)$-form.
We denote
\begin{align*}
	\theta_k\coloneqq& t_k\theta-(1-t_k)A\omega_{k},\quad\tau_k\coloneqq \Ric(\mu_{k})-\theta, \quad \beta_k\coloneqq \theta-\theta_k=(1-t_k)(\theta+A\omega_{k})\geq0,
\end{align*}
and
\begin{align*}
	\gamma_k\coloneqq\beta_k+\tau_k
	=\Ric(\mu_{k})-\theta_k.
\end{align*}
Let $M_k=M_{\eps_k,\mu_{\delta_k}}$ be the Mabuchi functional in $[\omega_{k}]$ twisting by $\gamma_k$, given by
\begin{align*}
	M_k(v)=\Ent_{\mu_k}(\omega_{k,v}^n/V_k)
	+c_kE_k(v)-E_{k,\theta_k}(v),
\end{align*}
and let $M_k^{\beta}=M_{\eps_k,\mu_{\delta_k}}^{\beta_k}$ be the Mabuchi functional in $[\omega_{k}]$ twisting by $\beta_k$, given by
\begin{align*}
	M_k^\beta(v)=\Ent_{\mu_k}(\omega_{k,v}^n/V_k)
	+\bar s_k^{\beta_k} E_k(v)-E_{k,\Ric(\mu_k)-\beta_k}(v),\quad
	\bar s_k^{\beta_k}\coloneqq
	\frac{n}{V_k}(c_1(Y)-[\beta_k])[\omega_k]^n.
\end{align*}
We also denote 
\begin{align*}
	\theta_t\coloneqq t\theta-(1-t)A\hat\omega, \qquad
	\beta_t\coloneqq\theta-\theta_t= (1-t)(\theta+A\hat\omega)\geq0,
\end{align*}
for some $A>0$.
  The $\beta_t$-twisted Mabuchi functional in $[\hat\omega]$ is given by
  \begin{align*}
  	M_t(v)=\Ent_{\mu}(\hat\omega_{v}^n/V)
	+c_tE_k(v)-E_{\theta_t}(v), \quad
	c_t\coloneqq\frac{n}{V}[\theta_t][\hat\omega]^{n-1}.
  \end{align*}
  
  \begin{proposition}
  \label{prop:Mab-compare}
  	Suppose $(u_k,F_k)$ is the solution of \eqref{eq:eps_twist_cscK} with a common entropy bound. Suppose $u\in\PSH_{\hat\omega}\cap L^\infty$ is limit of subsequence $u_k$, $\omega_k\rightarrow\hat\omega_u$ and $F_k\rightarrow F$ in $C^\infty_{\rm loc}(Y\setminus E)$.
  	Then for any $w\in\cH$,
  	\begin{align*}
  		M_t(u)\leq M_t(w).
  	\end{align*} 
  \end{proposition}
  
   On $Y\setminus E$, taking the limit for the first equation of \eqref{eq:eps_twist_cscK}, we have $\hat\omega_u^n=e^F\pi^*\mu_X$. Both measures have zero mass for $E$. Then the equation can be viewed as an equation of measures on $Y$. 
  Since $u,F$ is bounded, then $M_t(u)<\infty$.

\begin{lemma}
	\label{lem:converg-entropy}
	As the above notations, we have the convergence of entropy
	\begin{align*}
		\Ent_{\mu_k}(\omega_{k,w}^n/V_k)\rightarrow
		\Ent_{\mu}(\omega_{k,w}^n/V)
	\end{align*}
    and
	\begin{align*}
		\Ent_{\mu_k}(\omega_{k}'^n/V_k)\rightarrow
		\Ent_{\mu}(\hat\omega_{u}^n/V)
		=\int_Y Fe^F\mu.
	\end{align*}
	\end{lemma}
\begin{proof}
	We denote the density $f_{k,w}\coloneqq\frac{\omega_{k,w}^n}{V_k\omega_Y^n}$ and $f_w\coloneqq\frac{\hat\omega_w^n}{V\omega_Y^n}$. Then $f_{k,w}\rightarrow f_w$ smoothly and hence $f_{k,w}\log f_{k,w}\rightarrow f_w\log f_w$ uniformly (with $0\log 0=0$).
	Note that
	\begin{align*}
		\Ent_{\mu_k}(\omega_{k,w}^n/V_k)
		=\int_Y\log \frac{f_{k,w}}{f_k}f_{k,w}\omega_Y^n
		=\int_Yf_{k,w}\log f_{k,w}\omega_Y^n
		  -\int_Yf_{k,w}\log f_{k}\omega_Y^n,
	\end{align*}
	where $f_k=f_{\delta_k}$ the density of $\mu_{\delta_k}$, which decreasingly converges to the density of $\pi^*\mu_X$, denoted by $f_0$ as $\delta_k\rightarrow0$.
	It suffices to show the convergence of the second term.
	Note that
	\begin{align*}
		|\log f_k|\leq C+\sum_i|a_i||\log(r_i+\delta_k)|
		\leq C \sum_i|a_i||\log r_i|,
	\end{align*}
	and $\log f_k\rightarrow \log f_0$ pointwise on $Y\setminus E$.
	Since $\sum_i|a_i||\log r_i|$ is integrable (due to the klt condition), then the dominated convergence theorem yields
	\begin{align*}
		\int_Yf_{k,w}\log f_{k}\omega_Y^n
		\rightarrow \int_Yf_{w}\log f_{0}\omega_Y^n,
	\end{align*} 
	which proves the entropy convergence for $w$.
	
	\bigskip
	
	One has
	\begin{align*}
		\Ent_{\mu_k}(\omega_{k}'^n/V_k)
		=\int_Y F_ke^{F_k}f_k\omega_Y^n.
	\end{align*}
	Since $\|F_k\|_{L^\infty}\leq K$ (due to Theorem~\ref{thm:unif-C0}), then $|F_ke^{F_k}|\leq K$.
	For any compact $V\subset Y\setminus E$, one has
	\begin{align*}
		\int_V F_ke^{F_k}f_k\omega_Y^n\rightarrow
		\int_V Fe^{F}f_0\omega_Y^n.
	\end{align*}
	By \eqref{eq:Lp-int-mu-delta} and H\"older inequality, we have
	\begin{align*}
		\left|\int_{Y\setminus V}F_ke^{F_k}f_k\omega_Y^n \right|
		\leq K\|f_k\|_{L^{p_0}}\cdot\omega_Y^n(Y\setminus V)
		\leq  K\omega_Y^n(Y\setminus V)^{\frac{1}{p_0^*}}.
	\end{align*}
	Similarly,
	\begin{align*}
		\left|\int_{Y\setminus V}Fe^{F}f_0\omega_Y^n \right|
			\leq  K\omega_Y^n(Y\setminus V)^{\frac{1}{p_0^*}}.
    \end{align*}
	For any $\epsilon>0$, choose $V\Subset Y\setminus E$ such that $K\omega_Y^n(Y\setminus V)^{\frac{1}{p_0^*}}\leq \epsilon/3$. Then there exists $k_0>0$ such that for any $k\geq k_0$,
	\begin{align*}
		&\left|\Ent_{\mu_k}(\omega_{k}'^n/V_k)-\Ent_{\mu}(\hat\omega_{u}^n/V) \right| \\
		\leq& \left|\int_V F_ke^{F_k}f_k\omega_Y^n-
		\int_V Fe^{F}f_0\omega_Y^n \right|
		+\left|\int_{Y\setminus V}F_ke^{F_k}f_k\omega_Y^n \right|
		+\left|\int_{Y\setminus V}Fe^{F}f_0\omega_Y^n \right|
		\leq \frac{\epsilon}{3}+\frac{\epsilon}{3}+\frac{\epsilon}{3}=\epsilon.
	\end{align*}
	We finish the proof.
\end{proof}

\begin{lemma}
\label{lem:tau_k-formula}
	As the above notations, we have
	\begin{align*}
		\tau_k=-\sum_ia_i P_{i,k}
	\end{align*}
	where
	\begin{align*}
		P_{i,k}\coloneqq \ddc\log(r_i+\delta_k)+\Theta_i
		=\frac{\delta_k}{r_i+\delta_k}\Theta_i+\frac{\delta_k r_i}{(r_i+\delta_k)^2}\sqrt{-1}\partial\log r_i\wedge\bp\log r_i.
	\end{align*}
	Moreover, there exists $K>0$ such that $P_{i,k}\geq-K\omega_Y$. And $P_{i,k}\rightarrow0$ in $C_{\rm loc}^\infty(Y\setminus E)$.
\end{lemma}
\begin{proof}
	On $Y\setminus E$,
	\begin{align*}
		\Ric(\mu_k)-\theta=-\ddc (\log f_k-\log f_0)
		=&-\sum_i a_i\ddc\Big(\log(r_i+\delta_k)-\log r_i \Big) \\
		=&-\sum_i a_i\Big(\ddc\log(r_i+\delta_k)+\Theta_i \Big).
	\end{align*}
	Since $\tau_k$, $P_{i,k}$ are smooth forms. Then the above equality holds on $Y$.
	\medskip
	
	Since $|\Theta_i|\leq K\omega_Y$ and $0\leq\frac{\delta_k}{r_i+\delta_k}\leq1$, then $P_{i,k}\geq-K\omega_Y$.
	On a compact $V\Subset Y\setminus E$, there exists $c_V>0$ such that $r_i\geq c_K>0$ for all $i$. Thus
	\begin{align*}
		\left\|\frac{\delta_k}{r_i+\delta_k}\right\|_{C^m(V)}
		\leq C_{m,V}\delta_k\rightarrow0, \quad
		\left\|\frac{\delta_k r_i}{(r_i+\delta_k)^2}\right\|_{C^m(V)}
		\leq C_{m,V}\delta_k\rightarrow0,
	\end{align*}
	for any $m>0$, which yields that $P_{i,k}\rightarrow0$ in $C_{\rm loc}^\infty(Y\setminus E)$.
\end{proof}

\begin{lemma}
\label{lem:scalar-tr-L1->0}
For any $u_k$, then we have
\begin{enumerate}[$(\rm i)$]
	\item $\bar s_k^{\beta_k}-s(\omega_{k}')+\tr_{\omega_k'}\beta_k
	=(\bar s_k^{\beta_k}-c_k)-\tr_{\omega_k'}\tau_k$;
	\item $\|\bar s_k^{\beta_k}-s(\omega_k')+\tr_{\omega_k'}\beta_k\|_{L^1(\omega_k'^n)}\rightarrow0$, as $k\rightarrow\infty$;
	\item $|J_{k,\tau_k}(u_k)|\rightarrow0$ and $|J_{k,\tau_k}(w)|\rightarrow0$.
\end{enumerate}
\end{lemma}
\begin{proof}
	(i). By \eqref{eq:eps_twist_cscK}, one has
	\begin{align*}
		s(\omega_k')
		=\tr_{\omega_k'}\Ric(\omega_k')
		=\tr_{\omega_k'}(\Ric(\mu_k)-\ddc F)
		=&-\Delta_{\omega_k'}F+\tr_{\omega_k'}(\theta+\tau_k) \\
		=&~c_k+\tr_{\omega_k'}(-\theta_k+\theta+\tau_k) \\
		=&~c_k+\tr_{\omega_k'}(\beta_k+\tau_k).
	\end{align*}
	Thus we obtain (i).
	\medskip
	
	(ii). Note that
	\begin{align*}
		\bar s_k^{\beta_k}-c_k
		=\frac{n}{V_k}(c_1(Y)-[\beta_k]-[\theta_k])\cdot[\omega_k]^{n-1}
		=&~\frac{n}{V_k}(c_1(Y)-[\theta])\cdot[\omega_k]^{n-1} \\
		=&~=\frac{n}{V_k}\sum_i(a_i[E_i])\cdot[\omega_k]^{n-1}\rightarrow0.
	\end{align*} 
	By (i) and Lemma~\ref{lem:tau_k-formula}, it suffices to show
	\begin{align*}
		\int_Y|\tr_{\omega_k'}P_{i,k}|\omega_k'^n\rightarrow0.
	\end{align*}
	Since $[P_{i,k}]=[\Theta_i]=[E_i]$, then
	\begin{align*}
		\int_XP_{i,k}\wedge\omega_k'^{n-1}=[E_i][\omega_k]^{n-1}\rightarrow [E_i][\hat\omega]^{n-1}=0.
	\end{align*}
	Denote $P_{i,k}^\pm=\max(\pm P_{i,k},0)$.
	Since $P_{i,k}\geq -K\omega_Y$, then $P_{i,k}^-\leq K\omega_Y$.
	
	For any $\epsilon>0$, there exists a compact $V\Subset Y\setminus E$ such that $\omega_Y^n(Y\setminus V)\leq \epsilon$.
	By Lemma~\ref{lem:tau_k-formula}, one has
	\begin{align*}
		\int_VP_{i,k}^-\wedge\omega_k'^{n-1}
		\leq \int_V|P_{i,k}\wedge\omega_k'^{n-1}|
		\leq \|P_{i,k}\|_{C^0(V)}\int_V\omega_Y\wedge\omega_k'^{n-1}\rightarrow0
	\end{align*}
	and
	\begin{align*}
		\int_{Y\setminus V}P_{i,k}^-\wedge\omega_k'^{n-1}
		\leq K\int_{Y\setminus V}\omega_Y^n\leq K\epsilon.
	\end{align*}
	Thus, we have
	\begin{align*}
		\int_YP_{i,k}^-\wedge\omega_k'^{n-1}\rightarrow0.
	\end{align*}
	In addition, note that
	\begin{align*}
		\int_Y P_{i,k}^+\wedge \omega_k'^{n-1}
		=\int_Y P_{i,k}\wedge \omega_k'^{n-1}+\int_Y P_{i,k}^-\wedge \omega_k'^{n-1}\rightarrow0.
	\end{align*}
	Then we obtain
	\begin{align*}
		\int_Y|\tr_{\omega_k'}P_{i,k}|\omega_k'^n
		=n\int_Y P_{i,k}^+\wedge \omega_k'^{n-1}+n\int_Y P_{i,k}^-\wedge \omega_k'^{n-1}\rightarrow0,
	\end{align*}
	which proves (ii).
	
	 \medskip
	 
	 (iii). Since $\|u_k\|_{L^\infty}\leq K$, then $|E_k(u_k)|\leq K$.
	 One has
	 \begin{align*}
	 	\frac{n}{V_k}[\tau_k][\omega_k]^{n-1}
	 	=-\frac{n}{V_k}\sum_ka_i[E_i][\omega_k]^{n-1}\rightarrow0.
	 \end{align*}
	 By the similar calculation in (ii), we have
	 \begin{align*}
	 	|E_{k,\tau_k}(u_k)|\leq K\sum_{j=0}^{n-1}\int_Y\Big|\tau_k\wedge\omega_{k}'^{j}\wedge\omega_k^{n-1-j}\Big|\rightarrow0.
	 \end{align*}
	 Thus we obtain $|J_{k,\tau_k}(u_k)|\rightarrow0$. Similarly, $|J_{k,\tau_k}(w)|\rightarrow0$.
\end{proof}

  \begin{proof}[Proof of \textup{Proposition~\ref{prop:Mab-compare}}]

  Let $v_s$ be the $C^{1,1}$-geodesic in $[\omega_k]$ joining $v_0=u_k$ to $v_1=w$.
  By \cite[Theorem~1]{Dar17b} \cite[Section~2.2]{Bern15}, then
  \begin{align*}
  	|\dot v_0^+|\leq \|u_k-w\|_{L^\infty}\leq K,
  \end{align*}
  where $\dot v_0^+=\lim_{s\rightarrow0^+}\frac{v_s-v_0}{s}$.
  
  For any fixed $k$, let $\beta_{k,\lambda}=\beta_k+\lambda\omega_Y>0$ 
  By \cite[Lemma~3.3]{BDL20} and the convexity of the twisted Mabuchi functional (\cite{BB17a}, \cite{BDL17}), then
  \begin{align*}
  	M_k^{\beta_{k,\lambda}}(w)-M_k^{\beta_{k,\lambda}}(u_k)
  	\geq&~ \frac{1}{V_k}\int_Y\dot v_0^+(\bar s_k^{\beta_{k,\lambda}}-s(\omega_k')+\tr_{\omega_k'}\beta_{k,\lambda})\omega_k'^n \\
  	\geq&~-\|\dot v_0^+\|_{L^\infty}\frac{1}{V_k}
  	     \|\bar s_k^{\beta_{k,\lambda}}-s(\omega_k')+\tr_{\omega_k'}\beta_{k,\lambda}\|_{L^1(\omega_k'^n)} \\
  	\geq&~-K\|\bar s_k^{\beta_{k,\lambda}}-s(\omega_k')+\tr_{\omega_k'}\beta_{k,\lambda}\|_{L^1(\omega_k'^n)}.
  \end{align*}
  Let $\lambda\searrow0$, by Lemma~\ref{lem:scalar-tr-L1->0} (ii), we have
  \begin{align*}
  	M_k^{\beta}(w)-M_k^{\beta}(u_k)
  	\geq&~ -K\|\bar s_k^{\beta_k}-s(\omega_k')+\tr_{\omega_k'}\beta_k\|_{L^1(\omega_k'^n)}
  	=-o(1).
  \end{align*}
  Note that $M_k^\beta$ is linear with respect to the twisted form $\beta$.
  Then by Lemma~\ref{lem:scalar-tr-L1->0} (iii),
  \begin{align}
  	M_k(u_k)\leq M_k^\beta(u_k)-J_{k,\tau_k}(u_k)
  	\leq&~ M_k^\beta(w)+o(1)-J_{k,\tau_k}(u_k) \nonumber\\
  	\leq&~M_k(w)-J_{k,\tau_k}(u_k) +o(1)-J_{k,\tau_k}(u_k) \nonumber \\
  	\leq&~M_k(w)+o(1).
  	\label{eq:M_k(u_k)-leq-M_k(w)}
  \end{align}
  Since $u,w\in L^\infty$, then
  \begin{align*}
  	E_k(u_k)\rightarrow E(u) \quad\text{ and }\quad E_k(w)\rightarrow E(w).
  \end{align*} 
  Since $-K\omega_Y\leq\theta_k\leq K\omega_Y$ and $\theta_k\rightarrow\theta_t$ smoothly, then
  \begin{align*}
  	E_{k,\theta_k}(u_k)\rightarrow E_{\theta_t}(u) \quad\text{ and }\quad E_{k,\theta_k}(w)\rightarrow E_{\theta_t}(w).
  \end{align*}
  By Lemma~\ref{lem:converg-entropy}, we have
  \begin{align*}
  	M_k(u_k)\rightarrow M_t(u) \quad\text{ and }\quad
  	M_k(w)\rightarrow M_t(w). 
  \end{align*}
  	Combining with \eqref{eq:M_k(u_k)-leq-M_k(w)}, we finish the proof.
  \end{proof}

\subsection{Existence of solutions for perturbed equations}
In this subsection, we solve the perturbed system \eqref{eq:eps_twist_cscK} with a common entropy bound, in terms of Schauder fixed point theorem.
We always assume $G$-coercivity of Mabuchi functional $M$ in $[\omega]$.

\bigskip

We denote $K_{\min}\coloneqq\inf_{\cH_0}M>-\infty$.

\begin{lemma}
\label{lem:M(w_r)-leq-min+1}
	For any $r\in(0,1]$, there exists $w_r\in\cH_0$ such that
	\begin{align*}
		M(w_r)\leq K_{\min}+r, \qquad J_{\chi_0}(w_r)\leq D,
	\end{align*}
	where $\chi_0\coloneqq \theta+A\hat\omega$ and $D>0$ is independent of $r$.
\end{lemma}
\begin{proof}
	For any $r\in(0,1]$, one chooses $\varphi_r\in\cH_0$ such that $M(\varphi_r)\leq K_{\min}+r$.
	By \eqref{eq:coercive-M}, then
	\begin{align*}
		\d_{1,G}(0,\varphi_r)\leq \lambda_c^{-1}(K_{\min}+1+K_c).
	\end{align*}
	Then there exists $\sigma_r\in G$ such that $\d_1(0,  \varphi_r^{\sigma_r})\leq \lambda_c^{-1}(K_{\min}+1+K_c)+1$. We set $w_r\coloneqq  v_r^{\sigma_r}$. By $G$-invariance Lemma~\ref{lem:G-inv-Mabuchi}, $M(w_r)=M(v_r)\leq K_{\min}+r$.
	
	On the other hand, there exists a constant $D_{\chi_0}=D(\chi_0,n)$ such that for any $\varphi\in\cH_0$,
	\begin{align*}
		|J_{\chi_0}(\varphi)|\leq D_{\chi_0}\d_1(0,\varphi).
	\end{align*}
	By taking $D\coloneqq D_0\lambda_c^{-1}(K_{\min}+1+K_c)+D_0$, we finish the proof.
\end{proof}

\begin{lemma}
\label{lem:J_{chi-0}-lowerbound}
	For any $\varphi\in\PSH_\omega\cap L^\infty$ with $\sup \varphi=0$, then there exist $\constant{D}{constD:bound-L1}>0$ and $D_{\theta}>0$ such that
	\begin{align*}
		J_{\chi_0}(\varphi)\geq \left(\frac{A}{n}-D_\theta \right)(-E(\varphi))-D_1.
	\end{align*}
\end{lemma}
\begin{proof}
	For simplicity, we denote $E(j)\coloneqq\frac{1}{V}\int_Y(-\varphi)\hat\omega_\varphi^j\wedge\hat\omega^{n-j}$. Then $0\leq E(0)\leq\cdots\leq E(n)$.
	By $L^1$-compactness of qpsh functions, there exists $D_{\ref{constD:bound-L1}}>0$ such that
	\begin{align*}
		E(0)=\frac{1}{V}\int_Y(-\varphi)\hat\omega^{n}\leq D_{\ref{constD:bound-L1}}.
	\end{align*}
	Then one has
	\begin{align*}
		J_{\hat\omega}(\varphi)\geq-\frac{E(\varphi)}{n}-\frac{D_{\ref{constD:bound-L1}}}{n}.
	\end{align*}
	In addition, $J_\theta(\varphi)\geq-D_\theta(-E(\varphi))$.
	Then we obtain
	\begin{align*}
		J_{\chi_0}(\varphi)=J_\theta(\varphi)+AJ_{\hat\omega}(\varphi)
		\geq \left(\frac{A}{n}-D_\theta \right)(-E(\varphi))-D_1.
	\end{align*}
\end{proof}

\begin{proposition}
\label{prop:ent-gap}
	For any $t\in[0,1)$,  suppose $u\in\PSH_{\hat\omega}\cap L^\infty$, which is limit of subsequence $u_k$ of solution of \eqref{eq:eps_twist_cscK} for $t_k\rightarrow t$ and $\Ent_{\mu_k}(\omega_{k,u_k}^n/V_k) \leq K_{ent}$, then $\Ent_{\mu}(\hat\omega_{u}^n/V)\leq H_*$  for some constant $H_*>0$ independent of $t$ and $K_{ent}$.
\end{proposition} 
\begin{proof}
	Write $r=1-t>0$. Since $M_t=M+(1-t)J_{\chi_0}$, by Proposition~\ref{prop:Mab-compare}, then
	\begin{align*}
		M(u)+rJ_{\chi_0}(u)=M_t(u)\leq M_t(w_r)=M(w_r)+rJ_{\chi_0}(w_r).
	\end{align*}
	By Lemma~\ref{lem:M(w_r)-leq-min+1}, we have
	\begin{align*}
		K_{\min}+rJ_{\chi_0}(u)
		\leq M(w_r)+rJ_{\chi_0}(w_r)
		\leq K_{\min}+r+rD,
	\end{align*}
	which yields that $J_{\chi_0}(u)\leq D+1$, and hence by choosing $A>nD_\theta$ and Lemma~\ref{lem:J_{chi-0}-lowerbound},
	\begin{align*}
		-E(u)\leq \frac{1+D+D_{\ref{constD:bound-L1}}}{\frac{A}{n}-D_\theta}
		\eqqcolon \constant{D}{constD:unif-bound-(-E(u))}.
	\end{align*}
	 
	 On the other hand, by Proposition~\ref{prop:Mab-compare} with $w=0$, one has
	 \begin{align*}
	 	M(u)+rJ_{\chi_0}(u)\leq M(0).
	 \end{align*}
	 Again by Proposition~\ref{lem:J_{chi-0}-lowerbound}, then
	 \begin{align*}
	 	M(u)\leq M(0)-rJ_{\chi_0}(u)
	 	\leq M(0)+D_{\ref{constD:bound-L1}}.
	 \end{align*}
	 By \eqref{eq:Ric(mu_X)-bound}, then $|E_\theta(u)|\leq nK_0(-E(u))$. Thus, we obtain
	 \begin{align*}
	 	\Ent_{\mu}(\hat\omega_{u}^n/V)
	 	=M(u)-\bar s E(u)+E_\theta(u)
	 	\leq M(0)+D_{\ref{constD:bound-L1}}+(|\bar s|+nK_0)D_{\ref{constD:unif-bound-(-E(u))}}
	 	\eqqcolon H_*
	 \end{align*}
\end{proof}

\begin{corollary}
\label{cor:eps-threshold}
	For any $0<r<1$, there exists a constant $\eps(r)>0$ such that for all $\eps\in(0,\eps(r))$ and $t\in [0,1-r]$, the perturbed equation \eqref{eq:eps_twist_cscK} admits no solution with entropy in $[H,H+1]$, for any $H>H_*$.
\end{corollary}
\begin{proof}
	If not, then there exists $\eps_k\searrow0$, $\delta_k=\eps_k$ and $t_k\in[0,1-r]$ such that \eqref{eq:eps_twist_cscK} admits a solution $u_k$ with $\Ent_{\mu_k}(\omega_{k}'^n/V_k)\in [H,H+1]$. Up to subsequence, By Theorem~\ref{thm:unif-C0}, \ref{thm:unif-higher-order} and taking the limit, there exists $u\in \PSH_\omega\cap L^\infty$ such that $u_k\rightarrow u$ in $C_{\rm loc}^\infty(Y\setminus E)$.
	By Lemma~\ref{lem:converg-entropy}, one has
	\begin{align*}
		\Ent_{\mu}(\hat\omega_{u}^n/V)
		=\lim_{k\rightarrow\infty}\Ent_{\mu_k}(\omega_{k}'^n/V_k)\in [H,H+1].
	\end{align*}
    Proposition~\ref{prop:ent-gap} yields that $\Ent_{\mu}(\hat\omega_{u}^n/V)\leq H_*<H$, a contradiction.
    We finish the proof.
\end{proof}

\bigskip

For the rest of the subsection, we will produce a solution at $t=1-r$ for any $\eps=\delta\in(0,\eps(r))$, $r>0$.

\medskip

Fix $\gamma\in (0,1)$ and $\eps\in(0,1)$, we denote the Banach space
\begin{align*}
	\cB\coloneqq\left\{\phi\in C^{2,\gamma}(Y):\int_Y\phi\omega_Y^n=0 \right\}.
\end{align*}
For any $\phi\in\cB$, by the Aubin--Yau theorem, Evans--Krylov and Schauder's estimates, there exists an unique solution $u_\phi\in \cH_{\omega_\eps}$ of the equation
\begin{align}
\label{eq:Aubin-Yau-equ}
	V_\eps^{-1}(\omega_\eps+\ddc u_\phi)^n=e^{\phi+Au_\phi}\mu_\eps, \quad F_\phi=\phi+Au_\phi, 
\end{align}
and for some $0<\gamma'<\gamma$,
\begin{align}
\label{eq:Aubin-Yau-estimate}
	\|u_\phi\|_{C^{4,\gamma'}(Y,\omega_Y)}
	\leq K=K(\|\phi\|_{C^{2,\gamma}(Y,\omega_Y)},\omega_\eps,A).
\end{align}
The map $\phi\mapsto u_\phi$ is continuous. Indeed, if $\phi_j\rightarrow\phi$ in $\cB$, by the uniform $C^{4,\gamma'}$-estimate, up to subsequence, $u_{\phi_j}\rightarrow u_\infty$ in $C^{4,\gamma''}$ for some $\gamma''<\gamma'$ and $u_\infty\in C^{4,\gamma''}$, which satisfying the equation \eqref{eq:Aubin-Yau-equ}. Then $u_\phi=u_\infty$ by the uniqueness.

\medskip

We denote $\Lambda_\eps\coloneqq \theta+A\omega_\eps$ and $c_{\Lambda_\eps,\eps}\coloneqq\frac{n}{V_\eps}\int_Y\Lambda_\eps\wedge\omega_\eps^{n-1}$.
For any $\phi\in\cB$, there exists an uniqueness solution in $\cB$, denoted by $T_\eps(\phi)$, of the following Poisson equation
\begin{align}
\label{eq:Poisson-equ}
	\Delta_{\omega_{\eps,u_\phi}}\Big(T_\eps(\phi)\Big)
	=\tr_{\omega_{\eps,u_\phi}}\Lambda_\eps- c_{\Lambda_\eps,\eps},
\end{align}
satisfying
\begin{align}
\label{eq:Schauder-Poisson}
	\|T_\eps(\phi)\|_{C^{4,\gamma'}}\leq K=K(\|\phi\|_{C^{2,\gamma}},\omega_\eps,\gamma').
\end{align}
By the uniqueness and estimates \eqref{eq:Aubin-Yau-estimate}, \eqref{eq:Schauder-Poisson}, the map $T_\eps:\cB\rightarrow\cB$ is continuous and maps bounded sets to relatively compact sets, due to the compact embedding $C^{4,\gamma'}\hookrightarrow C^{2,\gamma}$.

\begin{lemma}
\label{lem:solu-equiv}
	For any $t\in[0,1)$, then $\phi=tT_\eps(\phi)$ holds iff $(u_\phi-\sup u_\phi, F_\phi)$ is a solution of \eqref{eq:eps_twist_cscK} at $\eps=\delta,t$. Such solution is smooth.
\end{lemma}
\begin{proof}
	If $\phi=tT_\eps(\phi)$, by \eqref{eq:Poisson-equ}, one has
	\begin{align*}
		\Delta_{\omega_{\eps,u_\phi}}F_\phi
		=\Delta_{\omega_{\eps,u_\phi}}\phi +A\Delta_{\omega_{\eps,u_\phi}}u_\phi 
		=&~t(\tr_{\omega_{\eps,u_\phi}}(\theta+A\omega_\eps)-c_{\Lambda,\eps})+A(n-\tr_{\omega_{\eps,u_\phi}}\omega_\eps) \\
		=&\tr_{\omega_{\eps,u_\phi}}\theta_{\eps,t}-c_{\eps_t}.
	\end{align*}
	
	If $(u,F)$ is a solution of \eqref{eq:eps_twist_cscK} at $\eps=\delta,t$.
	We set $\phi=F-Au-b$, where $b\coloneqq\frac{\int_Y(F-Au)\omega_Y^n}{\int_Y\omega_Y^n}$. Then $\phi\in\cB$. We also $\tilde u\coloneqq u+b/A$. Then $\phi+A\tilde u=F$ and $\tilde u=u_\phi$ satisfies \eqref{eq:Aubin-Yau-equ} by the first equation of \eqref{eq:eps_twist_cscK}.
	Moreover, one has
	\begin{align*}
		\Delta_{\omega_{\eps,u_\phi}}\phi
		=\Delta_{\omega_{\eps,u}}F-A\Delta_{\omega_{\eps,u}}u
		=&~\tr_{\omega_{\eps,u}}(t\theta-(1-t)A\omega_\eps)-c_{\eps,t}-A(n-\tr_{\omega_{\eps,u}}\omega_\eps) \\
		=&~t\tr_{\omega_{\eps,u_\phi}}(\theta+A\omega_\eps)-tc_{A_\eps,\eps}.
	\end{align*}
	Thus we have $\phi=tT_\eps(\phi)$.
\end{proof}

We denote $H_\eps(\phi)$ the entropy of $F_\phi$ with respect to $\mu_\eps$, given by
\begin{align*}
	H_\eps(\phi)\coloneqq\int_Y F_\phi e^{F_\phi}\mu_\eps.
\end{align*}
In particular, for $\phi=0$, set $u=u_0-\sup u_0$. Since $\int_Ye^{F_0}\mu_\eps=1$, then
\begin{align*}
	F_0=Au-\log\int_Ye^{Au}\mu_\eps.
\end{align*}
By Skoda integrability of $\kappa\omega_Y$-psh, H\"older inequality and \eqref{eq:eps_twist_cscK}, then
\begin{align*}
	\int_Y(-u)\mu_a\leq K_{\ref{constK:int(-u)e^F}}'=K(\omega_Y,\|f_\delta\|_{L^{p_0}},p_0).
\end{align*}
Jensen's inequality yields
\begin{align*}
	\log\int_Y e^{Au}\mu_\eps\geq A\int_Yu\mu_\eps\geq -AK_{\ref{constK:int(-u)e^F}}',
\end{align*}
which gives that $\sup F_0=-\log\int_Ye^{Au}\mu_\eps\leq AK_{\ref{constK:int(-u)e^F}}'$.
\begin{align}
	|H_\eps(0)|=\left|\int_Y F_0 e^{F_0}\mu_\eps\right|\leq 
	 AK_{\ref{constK:int(-u)e^F}}'+\frac{1}{e}\eqqcolon H_0.
\end{align}

\begin{theorem}
\label{thm:exist-perturb-equ}
Suppose Mabuchi functional $M$ is $G$-coercive.
	For any $r>0$ and any $\eps\in(\eps,\eps(r))$, then there exists a solution of \eqref{eq:eps_twist_cscK} at $\eps=\delta$, $t=1-r$, with uniform entropy bound.
\end{theorem}
\begin{proof}
Choose $H>\max(H_*,H_0)$. For a given $r\in(0,1)$, one takes $\eps\in(0,\eps(r))$ as in Corollary~\ref{cor:eps-threshold}.

\begin{claim}
\label{cla:phi-estimate}
	Suppose $\phi=tT_\eps(\phi)$ with $0\leq t\leq 1-r$, satisftying $H_\eps(\phi)\leq H+1$. Then there exists a constant $R_{\eps,r}>0$ such that $\|\phi\|_{C^{2,\gamma}}\leq R_{\eps,r}$.
\end{claim}
By Lemma~\ref{lem:solu-equiv}, then $\phi=F-Au-b$ for a solution $(u,F)$ of \eqref{eq:eps_twist_cscK} with entropy bound $H_\eps(\phi)=\int_YFe^F\mu_\eps\leq H+1$, where $b=\frac{\int_Y(F-Au)\omega_Y^n}{\int_Y\omega_Y^n}$. 
By Theorem~\ref{thm:unif-higher-order} (ii), one has
\begin{align*}
				\|u\|_{C^m(Y,\omega_Y)}+\|F\|_{C^m(Y,\omega_Y)}\leq K_{\eps,H,m}.
			\end{align*}
Then we obtain
\begin{align*}
	\|\phi\|_{C^{2,\gamma}}\leq C_\gamma(\|F-Au\|_{C^3}+|b|)
	\leq R_{\eps,r}.
\end{align*}

\medskip

Let $\eta_s$ be the continuous cutoff function, given by
\begin{align*}
	\eta_s(x)=\begin{cases}
		1, & \text{if }x\in(-\infty, s] \\
		-x+s+1, & \text{if }x\in[s,s+1] \\
		0, &\text{if }x\in[s+1,\infty).
	\end{cases}
\end{align*}
We define the map
\begin{align}
\label{eq;Schauder-map}
	S_{\eps,r}(\phi)\coloneqq(1-r)\eta_H\Big(H_\eps(\phi)\Big)\eta_{R_{a,r}}\Big(\|\phi\|_{C^{2,\gamma}}\Big) T_\eps(\phi).
\end{align}
Since $\eta_H$, $H_\eps(\cdot)$, $T_\eps(\cdot)$ is continuous, then the map $S_{\eps,r}$ is continuous.
\begin{claim}
	the image $S_{a,r}(\cB)$ is relatively compact.
\end{claim}
We denote $R\coloneqq R_{\eps,r}+1$ and $\cB_R\coloneqq\{\phi\in\cB:\|\phi\|_{C^{2,\gamma}}\leq R\}$. We first show that $V_R\coloneqq\bar{T_\eps(\cB_R)}$ is compact.
By \eqref{eq:Schauder-Poisson}, one has
\begin{align*}
	\sup_{\phi\in\cB_R}\|T_\eps(\phi)\|_{C^{4,\gamma'}}
	\leq K(\eps,R).
\end{align*}
For any sequence $\phi_j\in\cB_R$, by Arzel\`a--Ascoli theorem, then there exists a subsequence $\phi_{j_k}$ and $v\in C^{2,\gamma}$ such that 
\begin{align*}
	\|T_\eps(f_{j_k})-v\|_{C^{2,\gamma}}\rightarrow0,
\end{align*}
which gives that $\int_Yv\omega_Y^n=0$. Thus $V_R$ is compact.

We denote $t(\phi)\coloneqq (1-r)\eta_H\big(H_\eps(\phi)\big)\eta_{R_{a,r}}\big(\|\phi\|_{C^{2,\gamma}}\big) $. Then $0\leq t(\phi)\leq 1-r$.
If $\phi\notin\cB_R$, then $S_{\eps,r}(\phi)=0$. If $\phi\in\cB_R$, then $S_{\eps,r}(\phi)=t(\phi)T_\eps(\phi)$ with $T_\eps(\phi)\in V_R$.
Consequencely, we have
\begin{align*}
	S_{\eps,r}(\cB)\subset\{t f:t\in[0,1-r], f\in V_R\}.
\end{align*} 
Since $[0,1-r]\times V_R$ is compact, then $\bar{S_{\eps,r}(\cB)}$ is compact.

\bigskip

Let $L\coloneqq 1+\max_{f\in \bar{S_{\eps,r}(\cB)}}\|f\|_{C^{2,\gamma}}<\infty$. Then $S_{\eps,r}(B_L)\subset K\subset B_L$.
By Schauder fixed point theorem, there exists a fixed point $\phi\in B_L$, that is $\phi=t(\phi)T_\eps(\phi)$ with $t(\phi)\in[0,r]$.
By Lemma~\ref{lem:solu-equiv}, then  $(u_\phi-\sup u_\phi, F_\phi)$ is a solution of \eqref{eq:eps_twist_cscK} at $\eps=\delta,t$. 

If $H_\eps(\phi)\geq H+1$, then $t(\phi)=0$ and $\phi=0$, contradicting with $H_\eps(0)\leq H$.
By Corollary~\ref{cor:eps-threshold},
thus $H_\eps(\phi)<H$. By Claim~\ref{cla:phi-estimate}, then $\phi\in \cB_R$, which yields that $t(\phi)=1-r$.

Therefore, we obtain a solution of \eqref{eq:eps_twist_cscK} at $\eps,t=1-r$, with uniform entropy bound by $\max(H_*,H_0)$.
\end{proof}

\subsection{Proof of Theorem~\ref{thm:main}}
In this subsection, we finish the proof of Theorem~\ref{thm:main}, based on the uniform estimate given in Section~\ref{sec:Uniform esitmates for perturbed equations} and Theorem~\ref{thm:exist-perturb-equ}.

\begin{proof}[Proof of \textup{Theorem~\ref{thm:main}}]
Taking a sequence $r_k\searrow0$ and $\eps_k\in(0,\eps(r_k))$. 
	By the assumption and Theorem~\ref{thm:exist-perturb-equ},  
	there exists a solution $(u_{\eps_k},F_{\eps_k})$ of \eqref{eq:eps_twist_cscK} at $\eps_k=\delta_k$, $t_k=1-r_k$, satisfying
	\begin{align*}
		\Ent_{\mu_{\eps_k}}(\omega_{\eps_k,u_{\eps_k}}^n/V_{\eps_k})\leq H
	\end{align*}
	for some $H>0$.
	Then by Theorem~\ref{thm:unif-C0}, \ref{thm:unif-higher-order}, there exists $u\in\PSH_{\omega}\cap L^\infty\cap C^\infty(Y\setminus E)$ and $F\in L^\infty\cap C^\infty(Y\setminus E)$ such that (up to subsequence) $u_{\eps}\rightarrow u$, $F_{\eps_k}\rightarrow F$ in $C^\infty_{\rm loc}(Y\setminus E)$.
	Since $t_k\rightarrow1$, then on $Y\setminus E=X_{\reg}$, we have
	\begin{align*}
		V_\omega^{-1}(\omega+\ddc u)^n=e^F\mu_X, \qquad
		\Delta_{\omega_u}F=\tr_{\omega_u}\Ric(\mu_X)-\bar s,
	\end{align*}
	which yields that $s(u)=\tr_{\omega_u}\Ric(\omega_u)=\bar s$.
	By Proposition~\ref{prop:Mab-compare}, then $u$ is a minimizer of $M$.
	
	\medskip
	
	Finally, we prove the continuity of $u$. 
	For the equation $\hat\omega_u^n=e^F\pi^*\mu_X$ on $Y\setminus E$, both measures have zero mass for $E$. Then the equation can be viewed as an equation of measures on $Y$. 
	Since $F$ is bounded, then
	\begin{align*}
		\int_X (e^Ff_0)^{p_0}\omega_Y^n\leq K\|f_0\|_{L^{p_0}}^{p_0}\leq K.
	\end{align*}
    By \cite[Theorem~A]{EGZ09} (also \cite{CC25}), we obtain that $u\in C^0(X)$.
\end{proof}

\renewcommand\refname{Reference}
	\bibliographystyle{alpha}
	\bibliography{Reference_scscK.bib}

@article {AJL23,
    AUTHOR = {Apostolov, Vestislav and Jubert, Simon and Lahdili, Abdellah},
     TITLE = {Weighted {K}-stability and coercivity with applications to
              extremal {K}\"{a}hler and {S}asaki metrics},
   JOURNAL = {Geom. Topol.},
  FJOURNAL = {Geometry \& Topology},
    VOLUME = {27},
      YEAR = {2023},
    NUMBER = {8},
     PAGES = {3229--3302},
      ISSN = {1465-3060,1364-0380},
   MRCLASS = {32Q20 (14J45 32J27 53C25 53C55 58J60)},
  MRNUMBER = {4668097},
       DOI = {10.2140/gt.2023.27.3229},
       URL = {https://doi.org/10.2140/gt.2023.27.3229},
}

@article {BB17a,
    AUTHOR = {Berman, Robert J. and Berndtsson, Bo},
     TITLE = {Convexity of the {$K$}-energy on the space of {K}\"{a}hler
              metrics and uniqueness of extremal metrics},
   JOURNAL = {J. Amer. Math. Soc.},
  FJOURNAL = {Journal of the American Mathematical Society},
    VOLUME = {30},
      YEAR = {2017},
    NUMBER = {4},
     PAGES = {1165--1196},
      ISSN = {0894-0347,1088-6834},
   MRCLASS = {32Q15 (53C55)},
  MRNUMBER = {3671939},
MRREVIEWER = {Ruadha\'{\i}\ Dervan},
       DOI = {10.1090/jams/880},
       URL = {https://doi.org/10.1090/jams/880},
}

@article {BBEGZ19,
    AUTHOR = {Berman, Robert J. and Boucksom, Sebastien and Eyssidieux,
              Philippe and Guedj, Vincent and Zeriahi, Ahmed},
     TITLE = {K\"ahler-{E}instein metrics and the {K}\"ahler-{R}icci flow on
              log {F}ano varieties},
   JOURNAL = {J. Reine Angew. Math.},
  FJOURNAL = {Journal f\"ur die Reine und Angewandte Mathematik. [Crelle's
              Journal]},
    VOLUME = {751},
      YEAR = {2019},
     PAGES = {27--89},
      ISSN = {0075-4102,1435-5345},
   MRCLASS = {53C55 (53C25 53C44)},
  MRNUMBER = {3956691},
MRREVIEWER = {Kai\ Zheng},
       DOI = {10.1515/crelle-2016-0033},
       URL = {https://doi.org/10.1515/crelle-2016-0033},
}

@article {BDL17,
    AUTHOR = {Berman, Robert J. and Darvas, Tam\'{a}s and Lu, Chinh H.},
     TITLE = {Convexity of the extended {K}-energy and the large time
              behavior of the weak {C}alabi flow},
   JOURNAL = {Geom. Topol.},
  FJOURNAL = {Geometry \& Topology},
    VOLUME = {21},
      YEAR = {2017},
    NUMBER = {5},
     PAGES = {2945--2988},
      ISSN = {1465-3060,1364-0380},
   MRCLASS = {53C55 (32U05 32W20)},
  MRNUMBER = {3687111},
MRREVIEWER = {Chi\ Li},
       DOI = {10.2140/gt.2017.21.2945},
       URL = {https://doi.org/10.2140/gt.2017.21.2945},
}

@article {BDL20,
    AUTHOR = {Berman, Robert J. and Darvas, Tam\'{a}s and Lu, Chinh H.},
     TITLE = {Regularity of weak minimizers of the {$K$}-energy and
              applications to properness and {$K$}-stability},
   JOURNAL = {Ann. Sci. \'{E}c. Norm. Sup\'{e}r. (4)},
  FJOURNAL = {Annales Scientifiques de l'\'{E}cole Normale Sup\'{e}rieure.
              Quatri\`eme S\'{e}rie},
    VOLUME = {53},
      YEAR = {2020},
    NUMBER = {2},
     PAGES = {267--289},
      ISSN = {0012-9593,1873-2151},
   MRCLASS = {32Q26 (32Q15 53C55)},
  MRNUMBER = {4094559},
MRREVIEWER = {Kai\ Zheng},
       DOI = {10.24033/asens.2422},
       URL = {https://doi.org/10.24033/asens.2422},
}

@article {Bern15,
    AUTHOR = {Berndtsson, Bo},
     TITLE = {A {B}runn--{M}inkowski type inequality for {F}ano manifolds
              and some uniqueness theorems in {K}\"{a}hler geometry},
   JOURNAL = {Invent. Math.},
  FJOURNAL = {Inventiones Mathematicae},
    VOLUME = {200},
      YEAR = {2015},
    NUMBER = {1},
     PAGES = {149--200},
      ISSN = {0020-9910,1432-1297},
       DOI = {10.1007/s00222-014-0532-1},
       URL = {https://doi.org/10.1007/s00222-014-0532-1},
}

@article {BEGZ10,
    AUTHOR = {Boucksom, S. and Eyssidieux, P. and Guedj,
              V. and Zeriahi, A.},
     TITLE = {Monge-{A}mp\`ere equations in big cohomology classes},
   JOURNAL = {Acta Math.},
  FJOURNAL = {Acta Mathematica},
    VOLUME = {205},
      YEAR = {2010},
    NUMBER = {2},
     PAGES = {199--262},
      ISSN = {0001-5962},
   MRCLASS = {32U40 (32Q20 32U15 32W20)},
  MRNUMBER = {2746347},
MRREVIEWER = {S\l awomir Dinew},
       DOI = {10.1007/s11511-010-0054-7},
       URL = {https://doi.org/10.1007/s11511-010-0054-7},
}

@article {BJT26,
    AUTHOR = {Boucksom, S\'ebastien and Jonsson, Mattias and Trusiani,
              Antonio},
     TITLE = {Weighted extremal {K}\"ahler metrics on resolutions of
              singularities},
   JOURNAL = {Comm. Pure Appl. Math.},
  FJOURNAL = {Communications on Pure and Applied Mathematics},
    VOLUME = {79},
      YEAR = {2026},
    NUMBER = {4},
     PAGES = {1073--1148},
      ISSN = {0010-3640,1097-0312},
   MRCLASS = {58E11 (32U05 53C55)},
  MRNUMBER = {5030607},
       DOI = {10.1002/cpa.70026},
       URL = {https://doi.org/10.1002/cpa.70026},
}

@article {CC21b,
	AUTHOR = {Chen, Xiuxiong and Cheng, Jingrui},
	TITLE = {On the constant scalar curvature {K}\"ahler metrics
	({II})---{E}xistence results},
	JOURNAL = {J. Amer. Math. Soc.},
	FJOURNAL = {Journal of the American Mathematical Society},
	VOLUME = {34},
	YEAR = {2021},
	NUMBER = {4},
	PAGES = {937--1009},
	ISSN = {0894-0347,1088-6834},
	MRCLASS = {53C55 (53C21)},
	MRNUMBER = {4301558},
	MRREVIEWER = {Kai\ Zheng},
	DOI = {10.1090/jams/966},
	URL = {https://doi.org/10.1090/jams/966},
}

@misc{CC17,
      title={On the constant scalar curvature {K}\"ahler metrics, apriori estimates}, 
      author={Xiuxiong Chen and Jingrui Cheng},
      year={2017},
      eprint={1712.06697},
      archivePrefix={arXiv},
      primaryClass={math.DG}
}

@article {CC25,
    AUTHOR = {Cho, Ye-Won Luke and Choi, Young-Jun},
     TITLE = {Continuity of solutions to complex {M}onge-{A}mp\`ere
              equations on compact {K}\"ahler spaces},
   JOURNAL = {Math. Ann.},
  FJOURNAL = {Mathematische Annalen},
    VOLUME = {393},
      YEAR = {2025},
    NUMBER = {1},
     PAGES = {807--830},
      ISSN = {0025-5831,1432-1807},
   MRCLASS = {32W20 (14J17 32Q20 32U20)},
  MRNUMBER = {4966573},
MRREVIEWER = {Yaxiong\ Liu},
       DOI = {10.1007/s00208-025-03268-6},
       URL = {https://doi-org.proxy-um.researchport.umd.edu/10.1007/s00208-025-03268-6},
}

@article{CDL26,
 author = {Jianchun Chu and Tam{\'a}s Darvas and Yaxiong Liu},
 title = {The inverse {Monge}--{Amp{\`e}re} flow on klt {Fano} varieties},
 year = {2026},
 howpublished = {Preprint, {arXiv}:2610.36414 [math.{CV}] (2026)},
 url = {https://arxiv.org/abs/2610.36414},
 journal = {arXiv:2610.36414}
}

@article {Dar17a,
    AUTHOR = {Darvas, Tam\'as},
     TITLE = {Metric geometry of normal {K}\"ahler spaces, energy
              properness, and existence of canonical metrics},
   JOURNAL = {Int. Math. Res. Not. IMRN},
  FJOURNAL = {International Mathematics Research Notices. IMRN},
      YEAR = {2017},
    NUMBER = {22},
     PAGES = {6752--6777},
      ISSN = {1073-7928,1687-0247},
   MRCLASS = {32Q15 (32Q25 32W20 53C55)},
  MRNUMBER = {3737320},
MRREVIEWER = {S\l awomir\ Ko\l odziej},
       DOI = {10.1093/imrn/rnw203},
       URL = {https://doi.org/10.1093/imrn/rnw203},
}

@article {Dar17b,
    AUTHOR = {Darvas, Tam\'as},
     TITLE = {Weak geodesic rays in the space of {K}\"ahler potentials and
              the class {$\mathcal{E}(X,\omega)$}},
   JOURNAL = {J. Inst. Math. Jussieu},
  FJOURNAL = {Journal of the Institute of Mathematics of Jussieu. JIMJ.
              Journal de l'Institut de Math\'ematiques de Jussieu},
    VOLUME = {16},
      YEAR = {2017},
    NUMBER = {4},
     PAGES = {837--858},
      ISSN = {1474-7480,1475-3030},
   MRCLASS = {53C55 (32U05 32W20)},
  MRNUMBER = {3680345},
MRREVIEWER = {Matthew\ B.\ Stenzel},
       DOI = {10.1017/S1474748015000316},
       URL = {https://doi.org/10.1017/S1474748015000316},
}

@article {DR17,
    AUTHOR = {Darvas, Tam\'{a}s and Rubinstein, Yanir A.},
     TITLE = {Tian's properness conjectures and {F}insler geometry of the
              space of {K}\"{a}hler metrics},
   JOURNAL = {J. Amer. Math. Soc.},
  FJOURNAL = {Journal of the American Mathematical Society},
    VOLUME = {30},
      YEAR = {2017},
    NUMBER = {2},
     PAGES = {347--387},
      ISSN = {0894-0347,1088-6834},
   MRCLASS = {32Q20 (14J45 32U05 32W20 58B20 58E11)},
  MRNUMBER = {3600039},
MRREVIEWER = {Eleonora\ Di Nezza},
       DOI = {10.1090/jams/873},
       URL = {https://doi.org/10.1090/jams/873},
}

@article {Dem85,
    AUTHOR = {Demailly, Jean-Pierre},
     TITLE = {Mesures de {M}onge-{A}mp\`ere et caract\'{e}risation g\'{e}om\'{e}trique des
              vari\'{e}t\'{e}s alg\'{e}briques affines},
   JOURNAL = {M\'{e}m. Soc. Math. France (N.S.)},
  FJOURNAL = {M\'{e}moires de la Soci\'{e}t\'{e} Math\'{e}matique de France. Nouvelle S\'{e}rie},
    NUMBER = {19},
      YEAR = {1985},
     PAGES = {124},
      ISSN = {0037-9484},
   MRCLASS = {32H35 (32C10 32F05)},
  MRNUMBER = {813252},
MRREVIEWER = {G. M. Khenkin},
}

@article {DNG18,
    AUTHOR = {Di Nezza, Eleonora and Guedj, Vincent},
     TITLE = {Geometry and topology of the space of {K}\"ahler metrics on
              singular varieties},
   JOURNAL = {Compos. Math.},
  FJOURNAL = {Compositio Mathematica},
    VOLUME = {154},
      YEAR = {2018},
    NUMBER = {8},
     PAGES = {1593--1632},
      ISSN = {0010-437X,1570-5846},
   MRCLASS = {53C55 (32Q15 32W20 53C25)},
  MRNUMBER = {3830547},
MRREVIEWER = {Azniv\ Kasparian},
       DOI = {10.1112/s0010437x18007170},
       URL = {https://doi.org/10.1112/s0010437x18007170},
}

@article {DNGG23,
    AUTHOR = {Di Nezza, Eleonora and Guedj, Vincent and Guenancia, Henri},
     TITLE = {Families of singular {K}\"{a}hler-{E}instein metrics},
   JOURNAL = {J. Eur. Math. Soc. (JEMS)},
  FJOURNAL = {Journal of the European Mathematical Society (JEMS)},
    VOLUME = {25},
      YEAR = {2023},
    NUMBER = {7},
     PAGES = {2697--2762},
      ISSN = {1435-9855,1435-9863},
   MRCLASS = {32Q15 (14D06)},
  MRNUMBER = {4612100},
       DOI = {10.4171/jems/1249},
       URL = {https://doi.org/10.4171/jems/1249},
}

@article {DNJL26,
    AUTHOR = {Di Nezza, Eleonora and Jubert, Simon and Lahdili, Abdellah},
     TITLE = {Weighted csc{K} metric ({II}): {T}he continuity method},
   JOURNAL = {J. Reine Angew. Math.},
  FJOURNAL = {Journal f\"ur die Reine und Angewandte Mathematik. [Crelle's
              Journal]},
    VOLUME = {832},
      YEAR = {2026},
     PAGES = {255--296},
      ISSN = {0075-4102,1435-5345},
   MRCLASS = {32Q15 (53C55)},
  MRNUMBER = {5037728},
MRREVIEWER = {Charles\ Cifarelli},
       DOI = {10.1515/crelle-2026-0004},
       URL = {https://doi.org/10.1515/crelle-2026-0004},
}

@article {EGZ09,
    AUTHOR = {Eyssidieux, Philippe and Guedj, Vincent and Zeriahi, Ahmed},
     TITLE = {Singular {K}\"{a}hler-{E}instein metrics},
   JOURNAL = {J. Amer. Math. Soc.},
  FJOURNAL = {Journal of the American Mathematical Society},
    VOLUME = {22},
      YEAR = {2009},
    NUMBER = {3},
     PAGES = {607--639},
      ISSN = {0894-0347,1088-6834},
   MRCLASS = {32Q20 (31C10 32J27 32Q25 32W20)},
  MRNUMBER = {2505296},
MRREVIEWER = {Zhou\ Zhang},
       DOI = {10.1090/S0894-0347-09-00629-8},
       URL = {https://doi.org/10.1090/S0894-0347-09-00629-8},
}

@article{GSS26,
 author = {Bin Guo and Jian Song and Jacob Sturm},
 title = {Ricci entropy, {RCD} structures and {K}\"ahler spaces},
 year = {2026},
 howpublished = {Preprint, {arXiv}:2609.28986 [math.{DG}] (2026)},
 url = {https://arxiv.org/abs/2609.28986},
 journal = {arXiv:2609.28986}
}

@article {HLiu23,
    AUTHOR = {Han, Jiyuan and Liu, Yaxiong},
     TITLE = {On the existence of weighted-csc{K} metrics},
   JOURNAL = {Adv. Math.},
  FJOURNAL = {Advances in Mathematics},
    VOLUME = {463},
      YEAR = {2025},
     PAGES = {Paper No. 110125, 54},
      ISSN = {0001-8708,1090-2082},
   MRCLASS = {53C55 (32J27 53C21)},
  MRNUMBER = {4857568},
       DOI = {10.1016/j.aim.2025.110125},
       URL = {https://doi-org.proxy-um.researchport.umd.edu/10.1016/j.aim.2025.110125},
}

@article{HLiu25,
 author = {Jiyuan Han and Yaxiong Liu},
 title = {{$\mathbb{G}$}-uniform weighted {K}-stability for models on klt varieties},
 year = {2025},
 howpublished = {Preprint, {arXiv}:2506.18039 [math.{DG}] (2025)},
 url = {https://arxiv.org/abs/2506.18039},
 journal = {arXiv:2506.18039}
}

@article {Kol98,
	AUTHOR = {Ko{\l}odziej, S{\l}awomir},
	TITLE = {The complex {M}onge-{A}mp\`ere equation},
	JOURNAL = {Acta Math.},
	FJOURNAL = {Acta Mathematica},
	VOLUME = {180},
	YEAR = {1998},
	NUMBER = {1},
	PAGES = {69--117},
	ISSN = {0001-5962,1871-2509},
	MRCLASS = {32F07 (32C17 35J60)},
	MRNUMBER = {1618325},
	MRREVIEWER = {M.\ Klimek},
	DOI = {10.1007/BF02392879},
	URL = {https://doi.org/10.1007/BF02392879},
}

@article {Li22b,
    AUTHOR = {Li, Chi},
     TITLE = {Geodesic rays and stability in the csc{K} problem},
   JOURNAL = {Ann. Sci. \'Ec. Norm. Sup\'er. (4)},
  FJOURNAL = {Annales Scientifiques de l'\'Ecole Normale Sup\'erieure.
              Quatri\`eme S\'erie},
    VOLUME = {55},
      YEAR = {2022},
    NUMBER = {6},
     PAGES = {1529--1574},
      ISSN = {0012-9593,1873-2151},
   MRCLASS = {32Q26 (14E05 32Q15)},
  MRNUMBER = {4517682},
MRREVIEWER = {Ahmed\ Lesfari},
}

@article {LTW21,
    AUTHOR = {Li, Chi and Tian, Gang and Wang, Feng},
     TITLE = {On the {Y}au-{T}ian-{D}onaldson conjecture for singular {F}ano
              varieties},
   JOURNAL = {Comm. Pure Appl. Math.},
  FJOURNAL = {Communications on Pure and Applied Mathematics},
    VOLUME = {74},
      YEAR = {2021},
    NUMBER = {8},
     PAGES = {1748--1800},
      ISSN = {0010-3640,1097-0312},
   MRCLASS = {32Q26 (14J45 53C55)},
  MRNUMBER = {4275337},
MRREVIEWER = {Sergiy\ Koshkin},
       DOI = {10.1002/cpa.21936},
       URL = {https://doi.org/10.1002/cpa.21936},
}

@article {LTW22,
    AUTHOR = {Li, Chi and Tian, Gang and Wang, Feng},
     TITLE = {The uniform version of {Y}au-{T}ian-{D}onaldson conjecture for
              singular {F}ano varieties},
   JOURNAL = {Peking Math. J.},
  FJOURNAL = {Peking Mathematical Journal},
    VOLUME = {5},
      YEAR = {2022},
    NUMBER = {2},
     PAGES = {383--426},
      ISSN = {2096-6075,2524-7182},
   MRCLASS = {32Q20 (14J45 26E30 32P05 32Q26 53C55)},
  MRNUMBER = {4492658},
MRREVIEWER = {Yuchen\ Liu},
       DOI = {10.1007/s42543-021-00039-5},
       URL = {https://doi.org/10.1007/s42543-021-00039-5},
}

@article {LWX19,
    AUTHOR = {Li, Chi and Wang, Xiaowei and Xu, Chenyang},
     TITLE = {On the proper moduli spaces of smoothable
              {K}\"{a}hler-{E}instein {F}ano varieties},
   JOURNAL = {Duke Math. J.},
  FJOURNAL = {Duke Mathematical Journal},
    VOLUME = {168},
      YEAR = {2019},
    NUMBER = {8},
     PAGES = {1387--1459},
      ISSN = {0012-7094,1547-7398},
   MRCLASS = {14J45 (14D20 14J10 53C25 53C55)},
  MRNUMBER = {3959862},
MRREVIEWER = {Indranil\ Biswas},
       DOI = {10.1215/00127094-2018-0069},
       URL = {https://doi.org/10.1215/00127094-2018-0069},
}

@article {PT26,
    AUTHOR = {Pan, Chung-Ming and T\^o, Tat Dat},
     TITLE = {Weighted csc{K} metrics on {K}\"ahler varieties},
   JOURNAL = {J. Lond. Math. Soc. (2)},
  FJOURNAL = {Journal of the London Mathematical Society. Second Series},
    VOLUME = {113},
      YEAR = {2026},
    NUMBER = {2},
     PAGES = {Paper No. e70442, 47},
      ISSN = {0024-6107,1469-7750},
   MRCLASS = {53C55 (32Q20 32W20 35A23)},
  MRNUMBER = {5026855},
       DOI = {10.1112/jlms.70442},
       URL = {https://doi.org/10.1112/jlms.70442},
}

@article{PTT23,
  title={Singular csc{K} metrics on smoothable varieties},
  author={Pan, Chung-Ming and T{\^o}, Tat Dat and Trusiani, Antonio},
  journal={arXiv preprint arXiv:2312.13653},
  year={2023}
}

@article {SSY16,
    AUTHOR = {Spotti, Cristiano and Sun, Song and Yao, Chengjian},
     TITLE = {Existence and deformations of {K}\"{a}hler-{E}instein metrics
              on smoothable {$\Bbb{Q}$}-{F}ano varieties},
   JOURNAL = {Duke Math. J.},
  FJOURNAL = {Duke Mathematical Journal},
    VOLUME = {165},
      YEAR = {2016},
    NUMBER = {16},
     PAGES = {3043--3083},
      ISSN = {0012-7094,1547-7398},
   MRCLASS = {53C55 (14J10 14J45 32Q20 53C25)},
  MRNUMBER = {3566198},
MRREVIEWER = {G.\ K.\ Sankaran},
       DOI = {10.1215/00127094-3645330},
       URL = {https://doi.org/10.1215/00127094-3645330},
}

@book {Sze14,
    AUTHOR={Sz{\'e}kelyhidi, G\'{a}bor},
     TITLE = {An introduction to extremal {K}\"{a}hler metrics},
    SERIES = {Graduate Studies in Mathematics},
    VOLUME = {152},
 PUBLISHER = {American Mathematical Society, Providence, RI},
      YEAR = {2014},
     PAGES = {xvi+192},
      ISBN = {978-1-4704-1047-6},
   MRCLASS = {53C55 (14L24 32Q20 53C25)},
  MRNUMBER = {3186384},
MRREVIEWER = {Andrew\ Bucki},
       DOI = {10.1090/gsm/152},
       URL = {https://doi.org/10.1090/gsm/152},
}

@article {Sze25,
    AUTHOR = {Sz{\'e}kelyhidi, G\'{a}bor},
     TITLE = {Singular {K}\"ahler-{E}instein metrics and {RCD} spaces},
   JOURNAL = {Forum Math. Pi},
  FJOURNAL = {Forum of Mathematics. Pi},
    VOLUME = {13},
      YEAR = {2025},
     PAGES = {Paper No. e24, 33},
      ISSN = {2050-5086},
   MRCLASS = {53C55 (53C21 53C25)},
  MRNUMBER = {4976550},
MRREVIEWER = {Kai\ Zheng},
       DOI = {10.1017/fmp.2025.10015},
       URL = {https://doi.org/10.1017/fmp.2025.10015},
}

@article {Tian94,
    AUTHOR = {Tian, Gang},
     TITLE = {The {$K$}-energy on hypersurfaces and stability},
   JOURNAL = {Comm. Anal. Geom.},
  FJOURNAL = {Communications in Analysis and Geometry},
    VOLUME = {2},
      YEAR = {1994},
    NUMBER = {2},
     PAGES = {239--265},
      ISSN = {1019-8385,1944-9992},
   MRCLASS = {32G13 (32L07 53C07 53C25 58G26)},
  MRNUMBER = {1312688},
MRREVIEWER = {Kimio\ Miyajima},
       DOI = {10.4310/CAG.1994.v2.n2.a4},
       URL = {https://doi.org/10.4310/CAG.1994.v2.n2.a4},
}

@book {TG00,
    AUTHOR = {Tian, Gang},
     TITLE = {Canonical metrics in {K}\"{a}hler geometry},
    SERIES = {Lectures in Mathematics ETH Z\"{u}rich},
      NOTE = {Notes taken by Meike Akveld},
 PUBLISHER = {Birkh\"{a}user Verlag, Basel},
      YEAR = {2000},
     PAGES = {vi+101},
      ISBN = {3-7643-6194-8},
   MRCLASS = {32Q20 (53C25 58E11)},
  MRNUMBER = {1787650},
MRREVIEWER = {Azniv Kasparian},
       DOI = {10.1007/978-3-0348-8389-4},
       URL = {https://doi.org/10.1007/978-3-0348-8389-4},
}

@article {Yau78a,
    AUTHOR = {Yau, Shing Tung},
     TITLE = {On the {R}icci curvature of a compact {K}\"ahler manifold and
              the complex {M}onge-{A}mp\`ere equation. {I}},
   JOURNAL = {Comm. Pure Appl. Math.},
  FJOURNAL = {Communications on Pure and Applied Mathematics},
    VOLUME = {31},
      YEAR = {1978},
    NUMBER = {3},
     PAGES = {339--411},
      ISSN = {0010-3640,1097-0312},
   MRCLASS = {53C55 (32C10 35J60)},
  MRNUMBER = {480350},
MRREVIEWER = {Robert\ E.\ Greene},
       DOI = {10.1002/cpa.3160310304},
       URL = {https://doi.org/10.1002/cpa.3160310304},
}

@article {Yau78b,
    AUTHOR = {Yau, Shing Tung},
     TITLE = {A general {S}chwarz lemma for {K}\"ahler manifolds},
   JOURNAL = {Amer. J. Math.},
  FJOURNAL = {American Journal of Mathematics},
    VOLUME = {100},
      YEAR = {1978},
    NUMBER = {1},
     PAGES = {197--203},
      ISSN = {0002-9327,1080-6377},
   MRCLASS = {32H25 (53C55)},
  MRNUMBER = {486659},
MRREVIEWER = {Samuel\ I.\ Goldberg},
       DOI = {10.2307/2373880},
       URL = {https://doi.org/10.2307/2373880},
}

\end{document}